\documentclass[1p,authoryear]{elsarticle}

\usepackage[pdftex,colorlinks]{hyperref}
\usepackage{hypernat}
\usepackage{amsmath,amssymb,amsthm,verbatim,amscd}
\usepackage[usenames]{color}
\usepackage{mathrsfs,amsfonts,latexsym} 
\usepackage{psfrag,epsfig}
\usepackage{graphics}

\newcommand\NN{\mathbb N} 
 \newcommand\PP{\mathbb P}

 \newcommand{\EE}{\mathbb{E}}

\def\1{1\!{\mathrm l}}

\newcommand\rank{\operatorname{rank}}

\theoremstyle{plain} 
\newtheorem{thm}{Theorem}
\newtheorem{lem}[thm]{Lemma} 
\newtheorem{prop}[thm]{Proposition}
\newtheorem{cor}[thm]{Corollary}

\newtheorem{hyp}{Assumption}
\theoremstyle{definition} 
\newtheorem*{rem}{Remark}

\begin{document}

\begin{frontmatter}

\title{Parameter identifiability in a class of random graph mixture models}

\author[UAF]{Elizabeth S. Allman}
\ead{e.allman@alaska.edu}

\author[CM]{Catherine Matias\corref{cor1}}
\ead{catherine.matias@genopole.cnrs.fr}

\author[UAF]{John A. Rhodes}
\ead{j.rhodes@alaska.edu}

\cortext[cor1]{Corresponding author}

\address[UAF]{Department of Mathematics and Statistics, University of Alaska Fairbanks, PO Box 756660, Fairbanks, AK 99775, U.S.A}
\address[CM]{CNRS UMR 8071, Laboratoire Statistique et G\'enome, 523, place des Terrasses de l'Agora, 91 000 \'Evry, FRANCE}

\begin{abstract} 
 We prove identifiability of parameters for a broad class of random graph mixture models.  These models are characterized by a partition of the set of graph nodes into latent (unobservable) groups. The connectivities between nodes are independent random variables when conditioned on the groups of the nodes being connected. In the binary random graph case, in which edges are either present or absent, these models are known as stochastic blockmodels and have been widely used in the social sciences and, more recently, in biology. Their generalizations to weighted random graphs, either in parametric or non-parametric form, are also of interest in many areas. Despite a broad range of applications, the parameter identifiability issue for such models is involved, and previously has only been touched upon in the literature. We give here a thorough investigation of this problem. Our work also has consequences for parameter estimation. In particular, the estimation procedure proposed by Frank and Harary for binary affiliation models is revisited in this article.           
\end{abstract}

\begin{keyword} Identifiability  \sep mixture model \sep random graph \sep stochastic blockmodel 
\MSC[2008] 62E10 \sep 62F99
\end{keyword}

% 62E10 Characterization and structure theory
% 62F99 Parametric inference

\end{frontmatter}

%%%%%%%%%%%%%%

\section{Introduction}

In modern   statistical  analyses, data is often  structured  using networks.  Complex   networks   appear   across   many   fields   of science,  including biology  (metabolic networks,  transcriptional  regulatory  networks, protein-protein interaction networks), sociology (social networks of acquaintance, or other connections between individuals),  communications (the Internet), and others.

The literature contains many random graph models which incorporate a variety of characteristics of real-world  graphs  (such  as scale-free  or small-world properties).  We refer to \citet{Newman} and the references therein for an interesting introduction to networks.\\

One of the earliest and most studied random graph models was formulated by \citet{ER1}. In this setup, binary random graphs are modeled as a set of independent and identically distributed Bernoulli edge variables over a fixed set of nodes. The homogeneity of this model led to the introduction of mixture versions to better capture heterogeneity in data.  Stochastic blockmodels \citep{Daudin, FH82,Holland_etal_83,SN97} were introduced in various forms, primarily in the social sciences \citep{White_76} to study relational data, and more recently in biology \citep{Picard_BMC}. In this context, the nodes are partitioned into latent groups (blocks) characterizing the relations between nodes.  Blockmodelling thus refers to the particular structure of the adjacency matrix of the graph (\textit{i.e.}, the matrix containing edge indicators).  By ordering the nodes by the groups to which they belong,  this matrix exhibits a block pattern. Diagonal and off-diagonal blocks,  respectively, represent  intra-group and inter-group connections. In the special case where blocks exhibit the same behavior within their type (diagonal or off-diagonal), we obtain a model with an \emph{affiliation} structure \citep{FH82}. 

Although the literature from the social sciences has focused mostly on binary relations, there is a growing interest in weighted graphs \citep{PNAS_Barrat,Newman_weighted}. Mixture models have also been considered in the case of a finite number of possible relations \citep{NS01}, and more recently with continuous edge variables \citep{Matias_Ambroise, Mariadassou_Robin}. Some variations that we shall not discuss here include models with covariates \citep{Tallberg05}, mixed membership models \citep{Airoldi,Latouche09}, and models with continuous latent variables \citep{Daudin10,Handcock}. We also note that  \citet{Newman_Leicht} proposed another  version of a binary mixture model, slightly different from the stochastic blockmodel considered here.

Many different parameter estimation procedures have been proposed for these models, such as Bayesian methods \citep{NS01,SN97}, variational Expectation-Maximization (EM) procedures \citep{Daudin, Picard_BMC}, online classification EM methods \citep{Zanghi08,Zanghi09} and more recently, direct mixture model based approaches \citep{Matias_Ambroise}. Consistency of all these procedures relies strongly on the identifiability of the model parameters. However, the literature on these models has not addressed this question in any depth. The trivial label-swapping problem is often mentioned: it is  well known that the parameters may be recovered only up to permutations on the latent class labels. Whether this is the only issue preventing unique identification of parameters from the distribution, however, has never been investigated. Given the complex form of the model parameterization, this is not surprising, as any such analysis seems likely to be very involved. \\

In earlier work,  \citep[Theorem 7]{ECJ}, the authors made a first step towards an understanding of the parameter identifiability issue in binary random graph mixture models. While that article addressed a variety of models with latent variables, the present one focuses more specifically on random graph mixtures,  giving parameter identifiability results for a broad range of such models.  Moreover, part of our work sheds some new light on parameter estimation procedures.

\citet{ECJ} emphasized the usefulness of an algebraic theorem due to  \citet{Kruskal76,Kruskal77} \citep[see also][]{John_on_Kruskal} to establish identifiability results in various models whose common feature is the presence of latent groups and at least three conditionally independent variables. Here, we rather focus on the family of random graph mixture models and explore various techniques to establish their parameters' identifiability. Thus while the method developed by \citet{ECJ} is presented in Section~\ref{sec:scheme} and finds further use in several arguments, it is only one of several techniques we use. The issue at the core of Kruskal's result is the decomposition of a 3-way array as a sum of rank one tensors. While there exist approximate methods of performing this decomposition \citep[see, \textit{e.g.,}][]{Tomasi}, we mention that this approach seems poorly-suited to explicitly recover the parameters from the distribution, and thus to construct estimation procedures.

Some of our results focus on moment equations, as did  those of \citet{FH82}, in one of the earliest works on binary affiliation models. In particular, we revisit some of their claims. The method consists in looking at the distribution of $K_n$,  a complete set of edge variables over a set of $n$ nodes. A natural question is then: What is the minimal value of $n$ such that the complete distribution over all edge variables (a potentially infinite set) is characterized by the distribution of $K_n$? Despite this question's simplicity, we are far from having a complete answer to it. When looking at finite state distributions (\textit{e.g.}, for binary random graphs), the knowledge of the distribution of $K_n$ is equivalent to the knowledge of a certain set of moments of the distribution. Expressing the moments in terms of parameters gives a nonlinear polynomial system of equations, which one uses to identify parameters. The uniqueness of solutions to those systems, up to label swapping on parameters, is the issue at stake for identifiability. 

For random graphs with continuous edge weights given by a parametric family of distributions we shall see that the information contained in the model might be recovered from the distribution of $K_n$ for very small values of $n$. In this case, we rely on classical results on the identifiability of the parameters of a multivariate mixture due to \citet{Teicher67}. Note that the main difference between classical mixtures and random graph mixtures is the non-independence of the variates. 

In contrast to the approach based on Kruskal's Theorem, both the method utilizing moment equations and the one relying on multivariate mixtures lead to practical estimation procedures. These  are further developed by \citet{Matias_Ambroise}.\\

In \cite{ECJ}, a large role was played by the notion of \emph{generic identifiability}, by which every parameter except those lying on a proper algebraic subvariety, are identifiable. In other words, in a parametric setting, the non-identifiable parameters are included in a subset  whose dimension is strictly smaller than the dimension of the full parameter space. Thus with probability one with respect to the Lebesgue measure, every parameter is identifiable. This notion of generic identifiability is important for finite mixtures of multivariate Bernoulli distributions \citep{ECJ,Carreira,Gyllenberg} and also for hidden Markov models \citep{ECJ,Petrie69}.
Here, we stress that some of our identifiability results are generic, while others are strict.\\

Finally, we note that our focus throughout will be on undirected graph models. While many of our results may be generalized to directed graphs, one must pay careful attention to the models' parametrization in doing so. For instance, some of the results would become simpler if the connectivities from group $q$ to group $l$ differed from those from group $l$ to group $q$, as symmetry in a model can have a strong impact
on identifiability questions. However, such asymmetric models require an increase in the number of parameters which may be excessive for data analysis. \\

This paper is organized as follows. Section~\ref{sec:models} presents the various random graph mixture models: with either binary or, more generally, finite-state edges; both parametric and non-parametric models for edges with continuous weights; and the particular affiliation variant of these models. Section~\ref{sec:binary} gives parameter identifiability results for binary random graphs. Note that the  affiliation model has to be handled separately. Section~\ref{sec:weighted} takes up  weighted random graphs,  in both  parametric and non-parametric variants.  All the proofs are postponed to Section~\ref{sec:proofs}. In particular, Section~\ref{sec:scheme} is devoted to a brief presentation of  Kruskal's result and our use of it in the proofs of Theorems~\ref{thm:binary_Qge3_nonaff} and \ref{thm:discrete}.

\section{Notation and models}\label{sec:models}
We consider a probabilistic model on undirected and possibly weighted graphs as follows. Let $n$
be a fixed number of nodes, with $Z_1,\ldots,Z_n$ independent identically distributed (i.i.d.) random variables, taking values in
$\mathcal{Z}= \{1,\ldots,Q\}$ for some $Q\ge 2$.  These random variables represent the  $Q$ groups the
nodes are partitioned among, and are used to introduce heterogeneity in the model. With $\pi_q=\mathbb{P}(Z_i=q) \in (0,1)$, so $\sum_q \pi_q =1$, the vector $\boldsymbol \pi=(\pi_q)$ thus gives the priors on the groups. 
Let
$\{X_{ij}\}_{1\leq i<j\le n}$ be random edge variables taking values in a state space $ \mathcal{X}$.
Conditional on $Z_1,\ldots,Z_n$, we assume that the 
edge variables $\{X_{ij}\}_{1\leq i<j\le n} $ are independent, and that the conditional 
distribution of $X_{ij}$ depends only on $Z_i$ and $Z_j$, the groups containing its endpoints.

We are interested in random graphs of various types: For binary random graphs, where 
$\mathcal{X}=\{0,1\}$, an absent edge is represented by $0$ and a present one by $1$. Random graphs whose edges may be of finitely many types
are modeled with $\mathcal{X}=\{1,\ldots, \kappa\}$, or equivalently, $\{0,\ldots, \kappa-1\}$. More general weighted random graphs are obtained when  $\mathcal{X}=\mathbb{N}$ or $\mathbb{R}^s, s\ge 1$.\\

In the binary state case,
the distribution of $X_{ij}$ conditional on $Z_i,Z_j$ follows a Bernoulli distribution with parameter $p_{Z_iZ_j} = \mathbb{P}(X_{ij}=1| Z_i,Z_j)$.  As we consider only undirected graphs, we  implicitly assume  equality of the parameters $p_{ql}=p_{lq}$, for all $1\le q,l\le Q$. 

More generally, in the finite state case, with $\mathcal{X}=\{1,\ldots,\kappa\}$,
the vector $\mathbf{p}_{Z_iZ_j}=(p_{Z_iZ_j}(1),\ldots, p_{Z_iZ_j}(\kappa))$ contains the values $p_{Z_iZ_j}(k) = \mathbb{P}(X_{ij}=k| Z_i,Z_j)$, for $ 1\leq k\leq
\kappa$, with $\sum_k p_{Z_iZ_j}(k) =1$. We also implicitly assume  equality of the vectors $\mathbf{p}_{ql}=\mathbf{p}_{lq}$, for all $1\le q,l\le Q$.  We introduce this model primarily as a tool in the study of continuously weighted random graphs, though it might be useful for studying relationships between nodes of different types (colors), or of varying but discrete strengths (viewing the states as ordered). 
Note that a related model is described by \citet{NS01}, where the authors consider more general relation types (not necessarily edges, whether directed or not) occurring between a pair of nodes. \\

In the weighted random graph case, edges may be viewed as either absent ($X_{ij}=0$) or  present ($X_{ij}\ne 0$), with those present having a \emph{weight}, namely a non-zero value in $\mathcal{X}= \mathbb{N}$, $\mathbb{R}$, or $\mathbb{R}^s$. The distribution of  $X_{ij}$ conditional on $Z_i,Z_j$ may be assumed to have either a parametric or non-parametric form. More precisely, we assume that the distribution
of  $X_{ij}$ conditional on $Z_i,Z_j$ is the probability measure 
$\mu_{Z_i,Z_j}$ on $\mathcal{X}$ given by
\begin{equation*}
  \mu_{ql} = (1-p_{ql})\delta_0 +p_{ql}F_{ql} , \quad 1\le q,l\le Q, 
\end{equation*}
where $p_{ql}\in (0,1]$ is a sparsity parameter,  $\delta_0$ is the Dirac mass at zero and $F_{ql}$ is a probability measure on $\mathcal{X}$ with density $f_{ql}$ with respect to either the counting measure on $\mathbb{N}$ or the Lebesgue measure on $\mathbb{R}$ or $\mathbb{R}^s$. 
We also implicitly assume $\mu_{ql}=\mu_{lq}$, for all $1\le
q,l\le Q$. 
In the parametric case, we assume moreover that $F_{ql}=F(\cdot,\theta_{ql})$ and $f_{ql}=f(\cdot, \theta_{ql})$
where the parameter $\theta_{ql}$ belongs to $ \Theta \subset\mathbb{R}^p$.
In the non-parametric case we assume $F_{ql}$ is absolutely continuous.

We shall always assume that $F_{ql}$ has no point mass at zero, otherwise the sparsity parameter $p_{ql}$ cannot be identified from the mixture $\mu_{ql}$. For instance, when considering Poisson weights, $f_{ql}$ is the Poisson density truncated at zero, 
\begin{equation*}
  f_{ql}(k)= \frac{\theta_{ql}^k}{k !} (e^{\theta_{ql}}-1)^{-1}, \quad  k \ge 1 .
\end{equation*}

A particular instance of these models is the affiliation one, which assumes additionally
only two distributions of connections between the edges, one for 
intra-group connections and another for inter-group connections. Thus the
binary state case of the affiliation model assumes
\begin{equation*}
p_{ql}=
  \begin{cases}
    \alpha & \text{ if } q=l,\\
    \beta & \text{ if } q\neq l,
  \end{cases}
\quad\text{ for all } q,l\in \{1,\ldots,Q\}.
  \end{equation*}
The affiliation model in the continuous observations case is described similarly with $\mu_{ql} =\mu_{\text{in}} 1_{q=l}+\mu_{\text{out}} 1_{q\neq l}$, for all $1\le
q,l\le Q$. More precisely, in the continuous parametric case, for all $q,l\in  \{1,\ldots,Q\}$ we set 
\begin{equation*}
p_{ql}=
  \begin{cases}
    \alpha & \text{ if } q=l,\\
    \beta & \text{ if } q\neq l,
  \end{cases}
 \quad \text{ and } \quad 
\theta_{ql}=
  \begin{cases}
    \theta_{\text{in}} & \text{ if } q=l,\\
    \theta_{\text{out}} & \text{ if } q\neq l.
  \end{cases}
\end{equation*}

\medskip

For all these models, we consider restrictions of the model distribution by focusing on a subset of the nodes.
We denote by $K_n$ the complete set of $\binom {n}{2}$ edge variables associated to a subset of $n$ nodes. Note that the distribution of these variables is independent of the choice of which $n$ nodes one considers. Also, while this notation is motivated by that used in graph theory, where $K_n$ denotes the complete graph on $n$ nodes, we emphasize that here $K_n$ is a set of random variables, and we are making no statement as to whether these edges are present or absent in any realization of our model.

\section{Binary random graphs}\label{sec:binary}

We first focus on models with binary edge states, considering the more general case with arbitrary connectivity parameters, followed by  affiliation models.

\subsection{The non-affiliation case}
When $\mathcal{X}=\{0,1\}$, a first result on identifiability of parameters was obtained by \citet{ECJ} for the special case of $Q=2$ groups. For completeness, we recall the statement here. 

\begin{thm}\label{thm:previous} \cite[Theorem 7]{ECJ}. 
The parameters $\pi_1,\pi_2 = 1 - \pi_1, p_{11}, p_{12},p_{22}$ of the random graph mixture model with binary edge state variables and $Q=2$ groups are identifiable, up to label swapping, from the distribution of $K_{16}$  
provided that the connectivity parameters $\{p_{11},p_{12},p_{22}\}$ are  distinct. 

In particular, the result remains valid when the group proportions $\pi_q$ are fixed.
\end{thm}

Note the assumption that $p_{11}\ne p_{22}$ limits this theorem to the strict non-affiliation case.

The proof of this theorem is based on a clever application of an algebraic result, due to \citet{Kruskal76,Kruskal77} \citep[see also][]{John_on_Kruskal}, that deals with decompositions of 3-way arrays. While generalizing the proof to more than two groups requires substantially more effort, the basic method still applies.   Here we prove the following theorem.

\begin{thm}\label{thm:binary_Qge3_nonaff}
  The parameters $\pi_q$, $1\le q\le Q$, and $p_{ql}=\mathbb{P}(X_{ij}=1|Z_i=q,Z_j=l)$, $1\le q\le l\le Q$, of the random graph mixture model with binary edge state variables and $Q\ge 3$ groups
  are generically identifiable, up to label swapping, from the
  distribution of $K_{m^2}$, when
  \begin{equation*}
  \left\{
      \begin{array}{ll}
  m\ge Q-1 + (Q+2)^2/4 & \text{ if } Q \text{ is even},\\
  m\ge Q-1 +  (Q+1)(Q+3)/4 & \text{ if } Q \text{ is odd}.       
      \end{array}
\right.
  \end{equation*}
Moreover, the result remains valid when the group proportions $\pi_q$ are fixed.
\end{thm}

Note that the stated number of nodes ensuring that parameters are generically identifiable from the distribution of the edges may not be optimal. In particular, when $Q=2$, the proof of this theorem is still valid, yet it gives a minimal  number of $m^2= 25$ nodes. This is larger  than the bound $16$ obtained in Theorem~\ref{thm:previous}, and that number may itself not be optimal.

Also, while Theorem \ref{thm:previous} gives exact restrictions on parameters producing identifiability, Theorem \ref{thm:binary_Qge3_nonaff} is not explicit about
the generic conditions. However, for any fixed $Q$ the argument in our proof does yield a straightforward, though perhaps lengthy, means of checking whether a particular choice of parameters meets the conditions. Among these is a requirement that the $p_{ql}$ be distinct, so
the theorem does not apply to the affiliation model.

Moreover, a careful reading of the proof of the theorem shows that its generic aspect concerns only the part of the parameter space with the connectivities $p_{ql}$. This enables us to conclude that even when considering subsets defined by restriction of the group proportions $\pi_{q}$ (for instance assuming the group proportions are fixed, or equal), the result remains valid. 

\subsection{The affiliation model}\label{sec:binary_affiliation}

In the particular case of the affiliation model, we can obtain results from arguments based on moments of the distribution. For a small number of nodes, one may obtain explicit formulas for the moments in terms of model parameters. By analyzing the solutions to this nonlinear multivariate polynomial system of equations, one can address the question of parameter identifiability, as well as develop estimation procedures.

\subsubsection{Relying on the distribution of $K_3$.}
\citet{FH82} presented a method for estimation of the parameters of the binary affiliation model based only on the distribution of triplet cycles  $(X_{ij},X_{jk},X_{ki})$, $1\le i<j<k\le n$, of edge variables. From an identifiability perspective, this corresponds to identifying the parameters from the distribution of $K_3$. 
They suggest estimation of the parameters by solving the empirical moment equations. However, they omit discussing uniqueness of the solutions to these equations, even though this issue is a delicate one for nonlinear equations. 

In the following, we first explore the use of  the distribution of only $K_3$ to identify model parameters.  As a consequence, we exhibit a new estimation procedure for the parameters.  \\

The distribution of a triplet $(X_{ij},X_{jk},X_{ki})$, is expressible in terms of the indeterminates $\alpha,\beta$ and  $\pi_q$s. Let us denote by $s_2$ and $s_3$ the sums of the squares and cubes of the $\pi_q$s and, more generally, let
\begin{equation*}
s_k=\sum_{q=1}^Q \pi_q^k.  
\end{equation*}
Then one easily computes \citep[see also][]{FH82} the moment formulas
\begin{eqnarray}
  m_1&=&\mathbb{E}(X_{ij})= s_2\alpha+(1-s_2)\beta, \label{eq:m_1}\\
 m_2&=& \mathbb{E}(X_{ij}X_{ik})= s_3\alpha^2+2(s_2-s_3)\alpha\beta+(1-2s_2+s_3)\beta^2,  \label{eq:m_2}\\
m_3 &=&\mathbb{E}(X_{ij}X_{ik}X_{jk})= s_3\alpha^3+3(s_2-s_3)\alpha\beta^2+(1-3s_2+2s_3)\beta^3, \label{eq:m_3}
\end{eqnarray}
which completely characterize the distribution of $(X_{ij},$ $X_{jk},X_{ki})$.  

Note that in the important case of a uniform node distribution, where $\pi_q=1/Q$ for all $q$, we have
$s_k=Q^{1-k}$. This implies $s_3=s_2^2$, and hence
 $m_2=m_1^2$, so these equations reduce to two independent ones. As a consequence, the claim by \cite{FH82} that it is then possible to estimate the three unknowns $Q,\alpha,\beta$ relying only on these moment equations is not correct. 

Still, there are indeed several situations in which parameters are identifiable from these moments, as we next discuss. \\

With $Q=2$ latent groups and a possibly non-uniform group distribution, there are 3 independent parameters in the affiliation model.
In this case, the three moments above are enough to identify parameters.  To show this, we first construct certain polynomials with roots at the connectivity parameters.  Since the construction easily extends to larger $Q$, we give it more generally.

\begin{prop}\label{prop:alphapoly} Consider the random graph affiliation mixture model with $Q\ge2$ groups and binary edge state variables, on $Q+1$ nodes. Then
the parameter $\alpha$  is a real root of the degree $\binom{Q+1}2$ univariate polynomial 
$$U_Q(X)=\EE \left ( \prod_{1\le i<j\le Q+1} (X-X_{ij})\right ).$$
The polynomial
$$V_Q(X,Y)=\EE  \Big( \Big (X+(Q-1)Y-\sum_{1\le i\le Q    }X_{i(Q+1)}      \Big )\prod_{1\le i<j\le Q} (X-X_{ij})\Big )$$
of degree $\binom Q2 +1$ in $X,$ and degree 1 in $Y,$ vanishes at $(X,Y)=(\alpha,\beta)$. Moreover, the coefficient of
$Y$ in $V_Q(\alpha,Y)$ is non-zero precisely when $\alpha\ne \beta$.
\end{prop}

The utility of  these polynomials is that from the
distribution of $K_{Q+1}$, the  polynomial $U_Q$ allows one to recover
at most $\binom {Q+1}2$ candidate values  for $\alpha$, and then for each such value $V_Q$  allows   one  to   recover  a   unique   candidate  for
$\beta$. While some of  these candidates  could be ruled out  as not
lying  in $(0,1)$,  we do  not  know when this  leaves a  unique
$\alpha$ and  $\beta$ for $Q\ge 3$.  In the case of  $Q=2$ groups, however, we
prove that these polynomials uniquely identify the parameters.

\begin{thm}
  \label{thm:affil_Q=2} In the  random graph affiliation mixture model
  with $Q=2$ groups and binary edge state variables, 
the parameter $\alpha$ is the unique real root of the polynomial 
\begin{equation*}
U_2(X)=X^3-3m_1X^2+3m_2X-m_3 .  
\end{equation*}
Moreover, as soon as $\alpha \neq \beta$, the parameter $\beta$ is the unique real root of the
polynomial $V_2(\alpha,Y)$ where 
$$V_2(X,Y)=X^2+XY-3m_1X-m_1Y+2m_2.$$
\end{thm}

Once $\alpha$  and $\beta$ are  uniquely identified, we  may determine
from  equation \eqref{eq:m_1}  the value  of $s_2$  (again  using that
$\alpha\ne    \beta$),   and   hence    $\pi_1$,   $\pi_2$, up   to
permutation. This proves the following corollary.

\begin{cor}\label{cor:affil_Q=2}
The parameters $\{\pi_1,\pi_2 = 1 - \pi_1 \}$, up to label swapping, and $\alpha,\beta$ of the random graph affiliation mixture model with $Q=2$ groups and binary edge state variables are strictly identifiable from the distribution of $K_3$ provided $\alpha\neq \beta$. 
\end{cor}

\paragraph*{Identifiability of  $\alpha$ and $\beta$ when  $Q$ and the
  $\pi_q$s are known}  
  When the $\pi_q$s are known, \citet{FH82}  suggested   solving  any  two  of   the  three  empirical
counterparts   of   equations   \eqref{eq:m_1},   \eqref{eq:m_2}   and
\eqref{eq:m_3},  leading  to  three  different methods  of  estimating
$\alpha$  and  $\beta$.  However, numerical
experiments convinced us that two equations are in general not sufficient
to uniquely determine the parameters. In fact, it is not immediately clear that
even with the three moment
equations   (either  the   theoretical  ones   for  the   question  of
identification, or their  empirical counterparts for estimation) a unique solution is determined.
Below we give explicit formulas for the solution to the system, which in most cases are even rational, involving no extraction of roots.
These can thus be easily used to construct estimators.

\begin{thm}\label{thm:poly_ST}
If  $m_2\ne m_1^2$, then $\boldsymbol \pi$   is    non-uniform and
we can recover the parameters $\beta$ and $\alpha$ via the rational formulas
\begin{eqnarray*}
\beta &=&\frac{(s_3-s_2s_3)m_1^3   +(s_2^3-s_3)m_2m_1+(s_3s_2-s_2^3)m_3 }{(m_1^2-m_2)(2s_2^3-3s_3s_2+s_3)},\\ 
\alpha &=&\frac {m_1+(s_2-1)\beta}{s_2}. 
\end{eqnarray*}
If $m_2=m_1^2$, then $\boldsymbol \pi$ is uniform and we have 
\begin{equation*}
\beta =m_1+\left (  \frac {m_1^3-m_3}{Q-1} \right )^{1/3} \quad
\text{and} \quad 
\alpha =Q m_1+(1-Q)\beta.
\end{equation*}
\end{thm}
Implicit in this statement is the fact that denominators in the above formulas are non-zero. Note that the uniform group prior case formula is used for estimation by \cite{Matias_Ambroise}.

We immediately obtain the following corollary.

\begin{cor}\label{cor:poly_ST}
 For any fixed and known values of $\pi_q\in(0,1)$, $1\le q \le Q$, both parameters $\alpha,\beta$ of the random graph affiliation model with binary edge state variables are identifiable from the distribution of $K_3$. 
\end{cor}

The proofs of the previous statements lead to an interesting polynomial in the moments, whose vanishing detects the Erd\H{o}s-R\'enyi model, corresponding to a single node group.

\begin{prop}
The moments of a random graph affiliation model with binary edge state variables,  $Q$ node states,  and $\alpha\ne \beta$ satisfy
$$2m_1^3  - 3m_1m_2 +m_{3}=0$$
if, and only if, $Q=1$.
\end{prop}

This proposition follows from expressing the moments in terms of parameters to see that
$$2m_1^3  - 3m_1m_2 +m_3=(\alpha-\beta)^3(2s_2^3  - 3s_2s_3 +s_3),$$
together with the determination in the proof of Lemma~\ref{lem:aA} in Section \ref{sec:moment}  that $2s_2^3 - 3s_2s_3  +s_3\ne  0$  when  $\pi_q>0$  for  more than one group $q$.

\subsubsection{Relying on the distribution of $K_4$}

We next investigate parameter identifiability from the distribution of the edge variables over more than $3$ nodes, paying particular attention to the case of  $n=4$ nodes.

\paragraph*{Necessary conditions for identifiability of the $\pi_q$s, when $Q$ is known}
First, we establish that for an affiliation model, if the $\pi_q$s are unknown and are to be recovered from the distribution of $K_n$, then one must look at at least $n=Q$ nodes. Note that this applies not only to the binary edge state model, but to more general weighted edge models as well.
\begin{prop}\label{prop:bound_on_m}
  In order to identify, up to label swapping, the parameters $\{\pi_q\}_{1\le q \le Q}$ from
  an affiliation  random graph  mixture distribution on  $K_n$ (either
  binary or weighted), it is necessary that $n\ge Q$.  
\end{prop}

The condition in this lemma is in general not sufficient to identify the $\pi_q$. Indeed, the binary edge state affiliation model with $Q=3$ has $4$ parameters. However, the set of distributions over $K_3$ has dimension at most $3$ (according to equations \eqref{eq:m_1},\eqref{eq:m_2} and \eqref{eq:m_3}), which is not sufficient to identify the $4$ parameters.

\begin{figure}[h]
\begin{center}
\includegraphics[height=1.8cm]{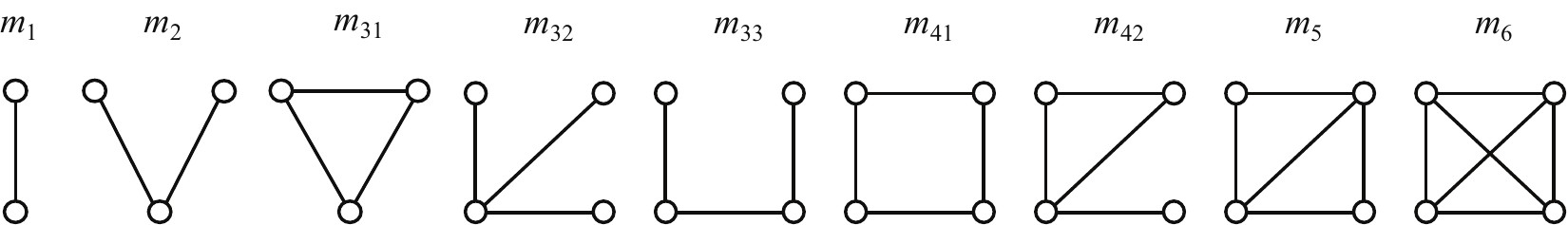} 
\end{center}
\caption{Correspondence between moments and motifs for $K_4$.}
\label{fig:moments}
\end{figure}

\paragraph*{Distribution on $K_4$}
The moment formulas describing the distribution of the affiliation random graph mixture model on $K_4$ are given in Table~\ref{fig:moment_K4}. Note that $m_{31}$ is the same as $m_3$ in the last subsection, and that we omit $\mathbb{E}(X_{12}X_{34})=(\mathbb{E}(X_{12}))^2$ since edge variables with no endpoints in common are independent. To facilitate understanding of the moments in the table, their corresponding induced motifs are shown in Figure~\ref{fig:moments}.

\begin{table}
    \begin{tabular}{|c|c|l|}
\hline
$m_1$ &$\mathbb{E}(X_{12})$&$ s_2\alpha +(1-s_2)\beta $\\
\hline
$ m_2 $ & $ \mathbb{E}(X_{12}X_{13})$ & $ s_3\alpha^2+2\alpha\beta(s_2-s_3)+(1-2s_2+s_3)\beta^2$\\
\hline
$ m_{31} $ & $ \mathbb{E}(X_{12}X_{13}X_{23})$ & $ s_3\alpha^3 +3(s_2-s_3)\alpha\beta^2+(1-3s_2+2s_3)\beta^3 $\\
\hline
$ m_{32} $ & $\mathbb{E}(X_{12}X_{13}X_{14}) $& $s_4\alpha^3 + 3(s_3-s_4)\alpha^2\beta+3(s_2-2s_3+s_4)\alpha\beta^2 $\\
& &$+(1-3s_2+3s_3-s_4)\beta^3$\\
\hline
$ m_{33} $ & $ \mathbb{E}(X_{12}X_{23}X_{34})$ & $s_4\alpha^3+(s_2^2+2s_3-3s_4)\alpha^2\beta+(3s_2-2s_2^2-4s_3+3s_4)\alpha\beta^2 $\\
&& $+(1-3s_2+s_2^2+2s_3-s_4)\beta^3$\\
\hline
 $m_{41} $ & $ \mathbb{E}(X_{12}X_{23} $ & $ s_4\alpha^4+2(s_2^2+2s_3-3s_4)\alpha^2\beta^2+4(s_2-s_2^2-2s_3+2s_4)\alpha\beta^3 $\\
&$X_{34}X_{41})$& $+(1-4s_2+2s_2^2+4s_3-3s_4)\beta^4$\\
\hline
$ m_{42} $ & $ \mathbb{E}(X_{12}X_{13}  $ & $ s_4\alpha^4+(s_3-s_4)\alpha^3\beta+(s_2^2+2s_3-3s_4)\alpha^2\beta^2 $\\
&$X_{14}X_{23})$& $ +(4s_2-2s_2^2-7s_3+5s_4)\alpha\beta^3 
 +(1-4s_2+s_2^2+4s_3-2s_4)\beta^4 $\\
\hline
 $m_5 $ & $ \mathbb{E}(X_{12}X_{23}X_{34} $ & $s_4\alpha^5+2(s_3-s_4)\alpha^3\beta^2+(2s_3-4s_4+2s_2^2)\alpha^2\beta^3 $\\ 
&$X_{41}X_{13})$ &$+(5s_2-4s_2^2-10s_3+9s_4)\alpha\beta^4
 +(1-5s_2+2s_2^2+6s_3-4s_4)\beta^5 $\\
\hline
$ m_6 $ &$\mathbb{E}(X_{12}X_{23}X_{34}$ & $s_4\alpha^6+4(s_3-s_4)\alpha^3\beta^3+3(s_2^2-s_4)\alpha^2\beta^4 $\\ 
&$X_{41}X_{13}X_{24})$&$ +6(s_2-s_2^2-2s_3+2s_4)\alpha\beta^5+(1-6s_2+8s_3-6s_4+3s_2^2)\beta^6 $\\
\hline
\end{tabular}
\caption{Moment formulas describing the distribution of $K_4$, the complete graph on $4$ nodes, for the binary affiliation model. }
\label{fig:moment_K4}
\end{table}

\smallskip

With $Q$ arbitrary, but a uniform prior on the nodes ($\pi_q=1/Q$, so $s_i=Q^{1-i}$), there are algebraic relationships between the moments on $K_4$, including
$$m_2=m_1^2,\ \  m_{32}=m_{33}=m_1^3,\ \  m_{42}=m_1m_{31},$$ and
more  complicated ones  that  can be  computed  using Gr\"obner  basis
methods  to eliminate $\alpha$, $\beta$, and $1/Q$ from the equations.  \citep[provide an excellent grounding on this  computational algebra.]{CLO} 
However,  the $3$ parameters $\alpha$, $\beta$, $Q$ of this affiliation model are, in fact, identifiable. Indeed such calculations show that the formulas for $m_1$, $m_{31}$, and $m_{41}$ alone imply the following.

\begin{prop}\label{prop:Qform}
The number of node  groups, $Q$, in a random graph affiliation model with binary edge state variables and uniform group priors can be identified from the moments  $m_1$, $m_{31}$, and $m_{41}$ by
$$Q= \frac{-m_{31}^4-m_{41}^3-3 m_{41}  m_1^8+3 m_{41}^2 m_1^4-6 m_1^6
  m_{31}^2+4 m_1^9 m_{31}+4 m_1^3 m_{31}^3}{(m_1^4-m_{41})^3}.$$
\end{prop}

Note that, replacing the moments with empirical estimators, this formula could be used for estimation of $Q$.
  
Of course once the formula in Proposition \ref{prop:Qform} is given, it can be most easily verified by expressing the moments in terms of parameters, and simplifying.
Note that  the denominator here does  not vanish, as may be seen in
two different ways: either by Lemma \ref{lem:4m} in Section \ref{sec:moment}, or
by checking that that $$m_{41}-m_1^4=(\alpha-\beta)^4\frac{(Q-1)}{Q^4}\ne 0.$$

Once $Q$ is identified by this formula, since we are assuming $\pi_q=1/Q$, Corollary \ref{cor:poly_ST} applies so that $\alpha$ and $\beta$ are identifiable as well. Thus we have shown the following.

\begin{cor} \label{cor:uniform_n=4}
The parameters  $\alpha$, $\beta$, and $Q$ of the random graph affiliation mixture model with binary edge state variables and uniform groups priors ($\pi_q=1/Q$) are identifiable from the distribution of $K_4$. 
\end{cor}

\section{Weighted random graphs}\label{sec:weighted}
\subsection{The parametric case}

 In the parametric case, where $F_{ql}$ has parametric form $F(\cdot, \theta_{ql})$, we can uniquely identify the connectivity parameters under very general conditions by considering the distribution of $K_3$ only. 
Indeed, each 
 triplet $(X_{ij},X_{ik},X_{jk})$ follows a mixture of  $Q^3$ distributions, each with three variates, comprising
 \begin{itemize}
 \item $Q$ terms of the form $\mu_{qq}(X_{ij})\mu_{qq}(X_{ik})\mu_{qq}(X_{jk}) $, each with prior $\pi_q^3$, where $1\le q\le Q$,
\item $3Q(Q-1)$ terms of the form $\mu_{qq}(X_{ij})\mu_{ql}(X_{ik})\mu_{ql}(X_{jk})$ (permuting $i,j$ and $k$), each with prior $\pi_{q}^2\pi_l$,  with distinct $ q,l\in \{1,2,\ldots,Q\}$,
\item $Q(Q-1)(Q-2)$ terms of the form $ \mu_{ql} (X_{ij}) \mu_{qm} (X_{ik}) \mu_{lm} (X_{jk})$,  each with prior $\pi_q\pi_l\pi_m$, with distinct $ q,l,m\in \{1,2,\ldots,Q\}$.
 \end{itemize}

By an old result due to \citet{Teicher67}, the identifiability of finite mixtures of some family of distributions is equivalent to
identifiability of finite mixtures of (multivariate) product distributions from this same family. In addition, identifiability of continuous
univariate parametric mixtures is generally well understood \citep{Teicher61,Teicher63}. Thus, we introduce the following assumptions.

\begin{hyp}\label{hyp:theta_neq}
  The $Q(Q+1)/2$  parameter values $\theta_{ql} $, $1\le q\le l\le Q$ are distinct.
\end{hyp}

\begin{hyp} \label{hyp:ident} The family of measures $\mathcal{M}=\{ F(\cdot,\theta) ~|~ \theta \in \Theta\}$ satisfies

\begin{enumerate}
\item[i)] all elements $F(\cdot,\theta)$ have no point mass at $0$,
\item[ii)] the parameters of finite mixtures of measures in $\mathcal{M}$ are identifiable, up to label swapping. In other words, for any integer $m\ge 1$, 
  \begin{equation*}
\text {if }  \quad 
\sum_{i=1}^{m} \alpha_i F(\cdot,\theta_i) =  \sum_{i=1}^{m} \alpha_i'  F(\cdot,\theta_i') 
\quad 
\text{ then } \quad \sum_{i=1}^{m} \alpha_i \delta_{\theta_i} = \sum_{i=1}^{m} \alpha_i' \delta_{\theta_i'} ,
  \end{equation*}
where $\delta_{\theta}$ denotes the Dirac mass at $\theta$.
\end{enumerate}
\end{hyp}

\begin{rem}
Note that most of the classical parametric families satisfy this assumption. In particular, the truncated Poisson, Gaussian and Laplace families $\{f(\cdot,\theta) , \theta \in \mathbb{R}^p\}$ satisfy Assumption~\ref{hyp:ident} \citep[see \emph{e.g.,}][]{MacLachlan_Peel,Teicher61,Teicher63}.
\end{rem}

\begin{thm}\label{thm:parametric}
  Under Assumptions~\ref{hyp:theta_neq} and~\ref{hyp:ident}, the parameters $\boldsymbol \pi$, $\theta_{ql}$, $p_{ql}$, $1\le q\le l\le Q$  of the parametric random graph mixture model with weighted edge variables are identifiable, up to label swapping, from the distribution of $K_3$.
\end{thm}

The previous result is not applicable to the parametric affiliation model, for which the set $\{\theta_{ql}, 1\le q\le l \le Q\}$ reduces to $\{\theta_{\text{in}}, \theta_{\text{out}}\}$, so Assumption~\ref{hyp:theta_neq} is violated. However, in this case a similar argument again yields a full identifiability result. As suggested by Proposition \ref{prop:bound_on_m}, we use $Q$ nodes to identify the group priors.

\begin{thm}\label{thm:parametric_affiliation} 
Under Assumption~\ref{hyp:ident}, the parameters $\alpha, \beta, \theta_{\text{in} },
  \theta_{\text{out}} $ of the parametric  affiliation  random graph mixture model with weighted edge variables are strictly  identifiable from the distribution of
$K_3$ provided $  \theta_{\text{in}}  \neq
\theta_{\text{out}} $. Once these have been identified, the  group priors $\boldsymbol \pi$  can further be
identified, up to label swapping, from the distribution of $K_Q$.
\end{thm}

A similar approach to that of this theorem has been  successfully used by \cite{Matias_Ambroise} to estimate the parameters of these models. They first estimated the sparsity parameters from the induced binary edge state model, but a procedure based on the preceding theorems would not require that these be distinct.\\

We turn next to models with a finite number, $\kappa$, of edge weights. Our primary reason for investigating such models is the role they play in our analysis of models with non-parametric conditional distributions of edge weights, in Section \ref{sec:nonpara}. Thus we limit our investigation to the single result we need there.

\begin{thm}\label{thm:discrete}
  The parameters of the random graph mixture model, with $\kappa$-state edge variables and $Q
  \ge  2$  latent groups, are identifiable, up to label
  swapping, from the
  distribution of $K_9$, provided $\kappa \ge
  \binom{Q+1}{2}$  and the $\kappa$-entry vectors  $\{\mathbf{p}_{ql}\}_{1\le q\le l\le Q}$ are linearly independent. 
\end{thm}

Note that the condition given here on the number of edge states is likely far from optimal. In
case $Q=2$ the condition requires at least $\kappa = 3$ edge states whereas we know from Theorem~\ref{thm:previous}  that the parameters are identifiable for this $Q$ with only $\kappa=2$ edge states.

\subsection{The non-parametric case}\label{sec:nonpara}

In the most general case of non-parametric distributions, our arguments for identifiability depend on binning the values of the edge variables into a finite set. We then apply Theorem \ref{thm:discrete} to this discretization, to obtain the following.

\begin{thm}\label{thm:nonparametric}
  The parameters $\{ \pi_q, \mu_{ql}=(1-p_{ql})\delta_0+p_{ql}F_{ql} : 1\le q,l\le Q\}$ of the random graph weighted non-parametric mixture model are identifiable, up to label swapping, from the distribution of $K_9$ provided the  measures $\mu_{ql} , 1\le q\le l \le Q $ are linearly independent.
\end{thm}

%%%%%%%%%%%%%%%%%%%%%%%%%%

\section{Proofs}\label{sec:proofs}

\subsection{Method of proofs based on Kruskal's theorem}\label{sec:scheme}

In this section we review Kruskal's theorem and describe our technique for employing
it in the proofs of Theorems~\ref{thm:binary_Qge3_nonaff} and \ref{thm:discrete}.

\paragraph*{Kruskal's result}

We first present Kruskal's result in a statistical context.  Consider
a latent random variable $V$ with state space $\{1,\ldots,r\}$ and
distribution given by the column vector $\mathbf v =(v_1,\ldots,
v_r)$. Assume that there are three observable random variables $U_j$
for $ j=1,2,3$, each with finite state space
$\{1,\ldots,\kappa_j\}$. The $U_j$s are moreover assumed to be
independent conditional on $V$.  Let $M_j$, $j=1,2,3$ be the
stochastic matrix of size $r\times \kappa_j$ whose $i$th row is
$\mathbf{m}_{i}^j = \PP(U_j=\cdot \mid V=i)$.  Then consider the
$\kappa_1\times\kappa_2\times \kappa_3$ tensor $[\mathbf
v;M_1,M_2,M_3]$ defined by
$$
[\mathbf v; M_1,M_2,M_3]=\sum_{i=1}^r v_i 
\mathbf m^1_i\otimes \mathbf m^2_i\otimes \mathbf m^3_i .
$$
Thus $[\mathbf v;M_1,M_2,M_3]$ is a $3$-dimensional array whose
$(s,t,u)$ element is
\begin{displaymath} 
[\mathbf v;M_1,M_2,M_3]_{s,t,u} =\sum_{i=1}^r v_i m^1_i(s) \,
  m^2_i(t) \, m^3_i(u) =\mathbb{P}(U_1=s,U_2=t,U_3=u),
\end{displaymath} for any $ 1\le s\le \kappa_1 , 1\le t\le \kappa_2 ,
1\le u\le \kappa_3 $.  
Note that $[\mathbf v;M_1,M_2,M_3]$ is left unchanged by
simultaneously permuting the rows of all the $M_j$ and the entries of
$\mathbf v$, as this corresponds to permuting the labels of the latent
classes.
Knowledge of the distribution of $(U_1,U_2,U_3)$ is equivalent to
knowledge of the tensor $[\mathbf v; M_1,M_2,M_3]$.\\

To state Kruskal's result, we need some algebraic terminology. For a
matrix $M$, the \emph{Kruskal rank} of $M$ will mean the largest
number $I$ such that \emph{every} set of $I$ rows of $M$ are
independent. Note that this concept would change if we replaced
``row'' by ``column,'' but we only use the row version in this
article. With the Kruskal rank of $M$ denoted by $\rank_K M$, we have
$$\rank_K M\le \rank M,$$
and equality of rank and Kruskal rank does not hold in general.
However, in the particular case when a matrix $M$ of size $p\times q$
has rank $p$, it also has Kruskal rank $p$.

The fundamental algebraic result of Kruskal is the following.

\begin{thm}\label{thm:kruskal}\citep{Kruskal76,Kruskal77}, \citep[see
  also][]{John_on_Kruskal} Let $I_j=rank_K M_j$. If
\begin{equation}
  \label{eq:kruskal_condition}
I_1+I_2+I_3\ge 2r+2,  
\end{equation}
then $[\mathbf v;M_1,M_2,M_3]$ uniquely determines $\mathbf v$ and the $M_j$, up to
simultaneous permutation of the rows. In other words, the set of
parameters $\{(\mathbf v, \mathbb{P}(U_j=\cdot \mid V))\}$ is uniquely
identified, up to label swapping, from the distribution of the random
variables $(U_1,U_2,U_3)$.
\end{thm}

Now, it will be useful to note that condition
\eqref{eq:kruskal_condition} holds for generic choices of the $M_j$,
provided the $\kappa_j$ are large enough to allow it. More precisely,
Kruskal's condition on the sum of Kruskal ranks can be expressed
through a Boolean combination of polynomial inequalities ($\neq$) involving matrix minors in
the parameters. If we show there is even a single choice of
parameters for which Kruskal's condition is satisfied, then the
algebraic variety of parameters for which it does \emph{not} hold is a proper
subvariety
(defined by negating the polynomial condition above, and so by a Boolean combination of equalities) of parameter space. As proper subvarieties are necessarily of 
Lebesgue measure zero, it follows that the Kruskal condition holds generically. \\

Our proof strategy for showing identifiability of certain random graph mixture models is to embed them in the
model we just described.  Applying Kruskal's result to the embedded
model, we derive partial identifiability results on the embedded model, and
then, using details of the embedding, relate these to the original model.

\paragraph*{Embedding the random graph mixture model into Kruskal's context}
Let $\kappa$ denote the cardinality of $\mathcal{X}$, in either the binary state case or the general finite state case.

To place the random graph mixture model in the context of Theorem
\ref{thm:kruskal}, we define a composite hidden variable and three composite observed
variables that reflect the conditional independence structure integral
to Kruskal's theorem.

For some $n$ (to be determined), let $V=(Z_1,Z_2,\dots,Z_n)$ be the latent random variable, with state
space $\{1,\ldots, Q\}^n$, which describes the state of all $n$ nodes
collectively, and denote by $\mathbf{v}$ the corresponding vector of
its probability distribution.  Note that the entries of $\mathbf{v}$
are of the form $\pi_1^{n_1}\cdots\pi_Q^{n_Q}$ with $n_q \ge 0$ and
$\sum n_q =n$.

The observed variables will correspond to three pairwise disjoint
subsets $G_1,G_2,G_3$ of the complete set of edges $K_n$. By choosing
the $G_i$ to have no edges in common, we ensure their conditional
independence.\\

The construction of the set of edges $G_i$ proceeds in two steps. We
begin by considering a small complete graph, and an associated matrix:
For a subset of $m$ nodes, we define a $Q^m \times
\kappa^{\binom m 2}$ matrix $A$, with rows indexed by assignments $\mathcal I\in
\{1,\ldots, Q\}^m$ of states to these $m$ nodes,
columns indexed by the state space of the complete set of $\binom m2$ edges between them, and
entries giving the probability of observing the specified states on all edges,
conditioned on the specified node states. In the case $\kappa=2$, it is helpful to note that each column index corresponds to a
different graph on the $m$ nodes, composed of those edges assigned state 1. For larger $\kappa$ one may similarly associate to a column index a $\kappa$-coloring of the edges of the complete graph. We therefore refer to a column index as a \emph{configuration}.

In the step we call the \emph{base case}, we exhibit a value of $m$ 
such that this matrix $A$ generically has full row rank. \\

Then, an \emph{extension step} builds on the base case, in order to
construct a larger set of $n$ nodes which will be used in the application of Theorem \ref{thm:kruskal}.
This is accomplished by means of \citep[][Lemma 16, and  subsequent remark]{ECJ} which we
paraphrase as follows.

\begin{lem} \label{lem:blockinduct} Suppose for the $Q$-node-state
  model, the number of nodes $m$ is such that 
  the $Q^m \times \kappa^{\binom{m}{2}}$ matrix $A$
  of probabilities of observing configurations of $K_m$ conditioned on node
  state assignments has rank $Q^m$. Then with $n=m^2$ there exist
  pairwise disjoint subsets $G_1,G_2,G_3$ of the complete set of edges $K_n$ such that
  for each $G_i$ the $Q^n\times \kappa^{|G_i|}$matrix $M_i$ of probabilities of observing
  configurations of $G_i$ conditioned on node state assignments has rank
  $Q^n$.
\end{lem}

In our applications here, we only determine that $A$  has full row rank generically. Hence the Lemma only allows
us to conclude that the $M_i$
have full row rank generically, and hence have Kruskal rank $Q^n$ generically. 

We also note (for use in the proof of Theorems \ref{thm:binary_Qge3_nonaff} and \ref{thm:discrete}) that in the construction of the Lemma,
each subset $G_j$ is the union of $m$
complete sets of edges each over $m$ different nodes, and thus
contains $m \binom m2$ edges.  In particular, if $m\ge 3$, then $G_i$ contains a complete graph on 3 nodes.

\paragraph*{Application of Kruskal's theorem to the embedded model and conclusion}
Next,
with
$\mathbf{v},M_1,M_2,M_3$ defined by the embedding given in the previous paragraphs, we apply Kruskal's Theorem (Theorem~\ref{thm:kruskal}) to the table $[\mathbf{v};M_1,M_2,M_3]$. Knowledge of the distribution of the random graph mixture
model over $n$ nodes implies knowledge of this $3$-dimensional table. By our construction of the $M_i$, condition
\eqref{eq:kruskal_condition} is satisfied since $3Q^n \geq 2Q^n +2$.  Thus
the vector $\mathbf{v}$ and the matrices $M_1,M_2,M_3$ are uniquely determined, up to simultaneous permutation of the rows.

With these embedded parameters in hand, it is still necessary to recover the initial parameters of the random graph mixture model: 
the group proportions and the connectivity vectors. As this requires a rather detailed argument, we leave its exposition for a specific application.

Finally, we note that by discretizing continuous variables, this
approach to establishing identifiability may also be used in the case of continuous connectivity
distributions.

\subsection{Proof of Theorem~\ref{thm:binary_Qge3_nonaff}}\label{sec:proof_binary_Qge3_nonaff}

This proof follows the strategy described in the previous section.  
We use the notation
$p_{ql}=\mathbb{P}(X_{ij}=1 \mid Z_i=q,Z_j=l)=1-\bar {p}_{ql}$.

\paragraph*{Base case} The initial step consists in finding a value of $m$
such that the matrix $A$ of size $Q^m\times 2^{\binom m 2}$
containing the probabilities of the configurations over these $m$ nodes, conditional on the hidden node states, 
generically has full row rank. 

The condition of having full row rank can be expressed as the
non-vanishing of at least one $Q^m \times Q^m$ minor of $A$. Composing
the map sending $\{p_{ql}\} \to A$ with this collection of minors gives polynomials
in the parameters of the model.  To see that these polynomials are
not identically zero, and thus are non-zero for generic parameters, it
is enough to exhibit a single choice
of the $\{p_{ql}\}$ for which the corresponding matrix $A$ has full row rank. 

With this in mind, we choose to consider $\{p_{ql}\}$ of the form
$p_{ql}=s_qs_l/(s_qs_l+t_qt_l)$, so
$\bar{p}_{ql}=t_qt_l/(s_qs_l+t_qt_l)$, with $s_i,t_j>0$ to be chosen
later.  However, since the property of having full row rank is unchanged
under non-zero rescaling of the rows of the matrix $A$, and all
entries of $A$ are monomials with total degree $\binom m 2$ in
$\{p_{ql},\bar{p}_{ql}\}$, we may simplify the entries of $A$ by 
removing denominators, and consider the matrix (also called
$A$) with entries in terms of $p_{ql}=s_q s_l$ and $\bar{p}_{ql}=t_q t_l$. 

The rows of $A$ are indexed by the composite node states
$\mathcal{I} \in \{1,\ldots,Q\}^m$, while its columns are indexed by
the edge configurations $\{0,1\}^{\binom m 2}$.  For any composite hidden state $\mathcal{I}
\in \{1,\ldots,Q\}^m$ and any vertex $v \in \{1,\ldots,m\}$, let
$\mathcal{I}(v) \in \{1,\ldots,Q\}$ denote the state of vertex $v$ in
the composite state $\mathcal{I}$.  With our particular choice of the
parameters $p_{ql}$, the $(\mathcal{I},(x_{ij})_{1\le i<j\le m})$-entry of $A$
is given by
\begin{equation*}
  \prod_{1\le v\le m} s_{\mathcal{I}(v)}^{d_v}
  t_{\mathcal{I}(v)}^{m-1-d_v} ,
\end{equation*}
where $d_v= \sum_{ w\neq v} x_{vw}$ is the degree of node $v$ in the
graph associated to the configuration $(x_{ij})_{1\le i<j\le m}$.  Note that the entries in a column of $A$
are now determined by the degree sequence $\mathbf{d}=(d_v)_{1\le
  v \le m}$ associated to the configuration. 
  
In general, there is a many-to-one correspondence of configurations to their 
degree sequences. (\textit{E.g.,} for $m=4$ nodes, the configuration with edges $(1,2)$ and $(3,4)$ in state 1, and that with
edges $(1,3)$ and $(2,4)$ in state 1, both have degree sequence $(1,1,1,1)$.)
Thus if $m>3$, there will be several identical columns in
$A$. 
For any degree sequence $\mathbf{d}=(d_v)_{1\le v \le m}$ arising
from an $m$-node graph, let $A_{\mathbf{d}}$ denote a
corresponding column of $A$.\\

Now, for each vertex $v\in \{1,\ldots,m\}$ and each $q\in
\{1,\ldots,Q\}$, introduce an indeterminate $U_{v,q}$ and a $Q^m$-entry row
vector $\mathbf{U}=(\prod_{1\le v\le m}
U_{v,\mathcal{I}(v)})_{\mathcal{I}\in \{1,\ldots,Q\}^m}$. For each
degree sequence $\mathbf{d}$, we have 
\begin{multline*}
  \mathbf{U}A_{\mathbf{d}}= \sum_{\mathcal{I}\in \{1,\ldots,Q\}^m}
  \prod_{1\le v\le m} s_{\mathcal{I}(v)}^{d_v}
  t_{\mathcal{I}(v)}^{m-1-d_v} U_{v,\mathcal{I}(v)} \\
= \prod_{1\le
    v\le m} \left(  s_{1}^{d_v}  t_{1}^{m-1-d_v}  U_{v,1} + \cdots +
s_{Q}^{d_v}  t_{Q}^{m-1-d_v}  U_{v,Q} \right).
\end{multline*}
To verify this,  notice
that each monomial $(s^{d_1}_{i_1} t^{m-1-d_1}_{i_1} U_{1,i_1}) \cdots (s^{d_m}_{i_m} t^{m-1-d_m}_{i_m} U_{m,i_m})$ obtained from multiplying out the product on the right corresponds to a choice of node states $i_v$ for nodes $v$, and hence a vector $\mathcal I = (i_1, \dots, i_m)$. Moreover, we obtain one such summand for each $\mathcal I$.

In order to prove that the matrix $A$ has full row rank, it is enough
to exhibit $Q^m$ independent columns of $A$.  Note, however, that
independence of a set of columns $\{A_\mathbf{d}\}$ is equivalent to
the independence of the corresponding set of polynomial functions $\{
\mathbf{U}A_{\mathbf{d}}\}$ in the indeterminates $\{U_{v,q}\}$.

Now for a set $\mathcal{D}$ of degree sequences, to prove that the
polynomials $\{\mathbf{U}A_{\mathbf{d}}\}_{\mathbf{d}\in \mathcal{D}}$
are independent, we assume that there exist scalars $a_{\mathbf{d}}$
such that
\begin{equation}\label{eq:null_sum}
  \sum_{\mathbf{d}\in \mathcal{D}} a_{\mathbf{d}}
  \mathbf{U}A_{\mathbf{d}} \equiv 0,
\end{equation}
and show that necessarily all $a_{\mathbf{d}}=0$. 
To this aim, we prove the following lemma.

\begin{lem}\label{lem:annihilate} Suppose $Q\le m$.
  Let $\mathcal{D}$ be a set of degree sequences such that for each
  node $v\in \{1,\ldots,m\}$, the set of degrees $\{d_v
  ~|~\mathbf{d}\in \mathcal{D}\}$ has cardinality at most $Q$. Then
  for generic values of $s_i,t_j$, for each $v$ and each $d^\star \in
  \{d_v~|~ \mathbf{d}\in \mathcal{D}\}$ there exist values of the
  indeterminates $\{U_{v,q}\}_{1\le q\le Q}$ that annihilate all the
  polynomials $\mathbf{U}A_{\mathbf{d}} $ for $\mathbf{d}\in
  \mathcal{D}$ except those for which $d_v=d^\star$.
\end{lem}

\begin{proof} Fix a node $v$ and let $\{d^1,\ldots,d^Q\}$ be any set
  of $Q$ distinct integers with
$$\{d_v~|~\mathbf{d}\in \mathcal{D}\} \subseteq\{d^1,\ldots,d^Q\}\subseteq\{0,1,\dots,m-1\}.$$

Let $M$ be the $Q\times
Q$ matrix with $i$th row $(s_1^{d^i}t_1^{m-1-d^i},
\ldots ,s_Q^{d^i}t_Q^{m-1-d^i})$.  Since all the integers $d^i$ are
different, the matrix $M$ has full row rank for generic choices of
$s_i,t_j$.   
(One way to see this is to consider a $m\times m$  
Vandermonde matrix, with $(k,l)$-entry $(u_l)^k$. Choosing distinct values of $u_l$ this has full rank, and thus the $Q\times m$ submatrix composed of rows with indices $\{d^i\}$ has rank $Q$. But then $Q$ of the columns can be chosen so that the $Q\times Q$
submatrix has full rank. Letting the $s_i$ be the values of $u_l$ in these columns, and $t_j=1$, gives one choice for which the matrix $M$ has full rank.)

Note  $d^\star=d^k$ for some $k$, and let
$\mathbf e_k$ be the $Q$-entry vector of all zeros except
for a $1$ in the $k$th position. Then for generic $s_i,t_j$, the equation
\begin{equation*}
  M (U_{v,1},\ldots,U_{v,Q})^{T} = \mathbf e_k\end{equation*}
admits a unique solution, one that corresponds to the above-mentioned choice of
indeterminates $\{U_{v,q}\}_{1\le q\le Q}$. 
\end{proof}

Now consider the following collection
\begin{multline*}
  \mathcal{D} =\Big\{ (d_1,\ldots,d_m) ~|~ d_v \in
  \{1,2, \ldots,Q\} \text{ for $v\le m-1$, and if }  \sum_{v=1}^{m-1} d_v \text{ is even}\\
  \text{ then } d_m \in \{0,2,4,\ldots,2Q-2\}, \text{ otherwise } d_m
  \in \{1,3,5,\ldots,2Q-1\} \Big\}.
\end{multline*}
Note that $\mathcal{D}$ has $Q^m$ elements and satisfies the
assumption of Lemma~\ref{lem:annihilate} on the number of different
values per coordinate. Moreover, if we establish, as we do below, that
its elements are realizable as degree sequences of graphs over $m$
nodes, then by choosing one column of $A$ associated to each degree sequence in $\mathcal{D}$,
we obtain a collection of $Q^m$ different columns of $A$. These
columns are independent since for each sequence $\mathbf{d}^\star\in
\mathcal{D}$ by Lemma~\ref{lem:annihilate} we can choose values of the
indeterminates $\{U_{v,q}\}_{1\le v\le m, 1\le q\le Q}$ such that all
polynomials $\mathbf{U}A_{\mathbf{d}}$ vanish, except
$\mathbf{U}A_{\mathbf{d}^\star}$, leading to $a_{\mathbf{d}^\star}=0$
in equation
\eqref{eq:null_sum}. \\

That each sequence $\mathbf{d}\in \mathcal{D}$ is realizable as a
degree sequence of a graph over $m$ nodes follows from a result of
 \citet{Erdos_Gallai} \citep[see also][Chapter 6,
Theorem 6]{Berge}.  Reordering the entries of $\mathbf d$ so that $d_1\ge d_2\ge
\ldots \ge d_m$, a necessary and sufficient condition for a sequence
to be realizable by such a graph is that for $1\le k\le m-1$,
\begin{equation}\label{eq:deg_condition}
 \sum_{v=1}^k d_v \le k(k-1)
  +\sum_{v=k+1}^m \min\{k,d_v\}.
\end{equation}
From the definition of $\mathbf d\in\mathcal D$, with coordinates
reordered, it is easy to see that for any $1\le k\le m-1$, we have
\begin{equation*}
  \sum_{v=1}^k d_v \le (k-1)Q +(2Q-1)\quad \text{and} \quad
\sum_{v=k+1}^m \min\{k,d_v\} \ge  m-k.
\end{equation*}
Thus, for \eqref{eq:deg_condition} to be satisfied, it is enough that
for any $1\le k\le m-1$, we have
\begin{equation*}
 - k^2 +(Q+2)k +Q-1\le m.
\end{equation*}
But for $m$ sufficently large
\begin{equation*}
  \max_{1\le k\le m-1} \{-k^2 +(Q+2)k \}= \left\{
    \begin{array}{ll}
\left( \frac{Q+2}{2}\right)^2 & \text{if } Q \text{ is even,}\\      
\frac{(Q+1)(Q+3)}{4} & \text{if } Q \text{ is odd.}     
    \end{array}
\right.
\end{equation*}
Thus, inequality \eqref{eq:deg_condition} is satisfied as soon as 
\begin{equation*}
  \left\{
      \begin{array}{ll}
  m\ge Q-1 + \big( \frac{Q+2}{2}\big)^2 & \text{ if } Q \text{ is even},\\
  m\ge Q-1 +  \frac{(Q+1)(Q+3)}{4} & \text{ if } Q \text{ is odd}.       
      \end{array}
\right.
  \end{equation*}
This concludes the proof of the base case.

The extension step explained in Section \ref{sec:scheme} then applies, so that with $n=m^2$, Kruskal's Theorem may be
applied to identify, up to simultaneous row permutation, $\mathbf v$, $M_1$, $M_2$, and $M_3$ as defined in that section.

\paragraph*{Conclusion}
The entries of $\mathbf{v}$ obtained via Kruskal's theorem
applied to the embedded model are of the form
$\pi_1^{n_1}\cdots\pi_Q^{n_Q}$ with $\sum n_q =n$, while the entries of the $M_i$ contain information on the $p_{ql}$. Although
the ordering of the rows of the $M_i$ is arbitrary, crucially we do know how the rows of $M_i$ are paired with the entries of $\mathbf v$.

By focusing on one of the matrices, say  $M_1$, and adding appropriate columns
to marginalize to a single edge variable (\emph{e.g.}, all columns for configurations with $x_{12}=1$), we recover the 
set of values $\{p_{ql}\}_{1\le q \le l\le Q}$, but without order.    However, if row $k$ of $M_1$ corresponds to the unknown
node states $\mathcal I$,
then performing 
such marginalizations for each of the $3$ edges of a complete graph $C$ on $3$ nodes contained in $G_1$ recovers the set
$$
R_k = \{ p_{ql} \mid  \text{ for some edge $(v,w)\in C$, $\{\mathcal I(v),\mathcal I(w)\}=\{q,l\}$ } \}.
$$
By considering the cardinalities of the sets $R_k$ in the generic case of all $p_{ql}$ distinct, we can now
determine individual parameters.   

Consider first those $k$ for which $R_k$ has one element.  There
are exactly $Q$ of these, arising from all $3$ nodes being in the same group.  
Thus for such $k$, $R_k=\{p_{qq}\}$ and $v_k=\pi_q^n$. Choosing an arbitrary labeling, we have determined
all $\pi_q$  and $p_{qq}$.

Next consider those $k$ for which the $R_k$ has two
elements.   These arise from $2$ nodes being in the same group, with the other node in a
different group, so $R_k=\{p_{qq},p_{ql}\}$ for some
$l \neq q$.  However, having already determined the $p_{qq}$
and since generically the $p_{ql}$ are distinct, we can
find exactly two such $k_1$ and $k_2$ of the form $R_{k_1} = \{ p_{qq}, p_{ql} \}$
and $R_{k_2} = \{ p_{ll}, p_{ql} \}$.  Thus, we can also determine $p_{ql}$
for $q \neq l$.

Finally, note that all generic aspects of this argument, in the base case and the requirement that the parameters $p_{ql}$ be distinct, concern only the $p_{ql}$. Thus if the group proportions $\pi_q$ are fixed to any specific values, the theorem remains valid. \qed

%%%%%%%

\subsection{Proofs relying on moment equations}\label{sec:moment}

\begin{proof}[Proof of Proposition~\ref{prop:alphapoly}] 
Focusing on $Q+1$ nodes, let
$Z=(Z_1,\ldots,Z_{Q+1})$ denote the composite node random variable, and  $z=(z_1,\ldots,z_{Q+1})$ any realization of $Z$.
Note that 
\begin{align*}U_Q(X)
&=\sum_{z\in \{1,\dots,Q\}^{Q+1}} \left (\prod_{1\le k\le Q+1} \pi_{z_k} \right ) \EE \left ( \prod_{1\le i<j\le Q+1} (X-X_{ij})~|~ Z=z\right )\\
&=\sum_{z\in \{1,\dots,Q\}^{Q+1}} \left (\prod_{1\le k\le Q+1} \pi_{z_k} \right )\prod_{1\le i<j\le Q+1} \left (X-\EE(X_{ij}~|~Z_i=z_i,Z_j=z_j) \right),
\end{align*}
since  conditioned   on  $Z=z$,  the   edge  variables
$X_{ij}$ are independent. Now since there are $Q+1$ nodes and only $Q$
groups, for each term in the sum there is some $z_i=z_j$. Since
$$ X- \EE (X_{ij}|Z_i=z_i=z_j=Z_j)=X-\alpha,$$
each term in the sum vanishes at $X=\alpha$, so $U_Q(\alpha)=0$.\\

Likewise,
\begin{align*}V_Q(X,Y)&=\sum_{z\in \{1,\dots,Q\}^{Q+1}} \left (\prod_{1\le k\le Q+1} \pi_{z_k} \right ) \times\\
&\qquad \EE \left (  \left. \left (X+(Q-1)Y-\sum_{1\le i\le Q    }X_{i(Q+1)}      \right )\prod_{1\le i<j\le Q} (X-X_{ij})~\right |~ Z=z\right ).
\end{align*}
But
\begin{multline*}
\EE    \left(  \left. \left    (X+(Q-1)Y-\sum_{1\le    i\le   Q    }X_{i(Q+1)}
\right )\prod_{1\le i<j\le Q} (X-X_{ij})~\right |~ Z=z\right )\\
=     \left (       X+(Q-1)Y-\sum_{1\le i\le  Q}
\EE\left ( X_{i(Q+1)}~|~Z_i=z_i,Z_{Q+1}=z_{Q+1}  \right )\right )\times \\
\prod_{1\le i<j\le Q} \left (X-\EE (X_{ij}~|~Z_i=z_i,Z_j=z_j) \right ).
\end{multline*}
Letting $X=\alpha$, one of the factors  $X-\EE (X_{ij}~|~ Z_i=z_i,Z_j=z_j)$ will vanish for any $z$ except possibly those with the $z_i$, $1\le i\le Q$, distinct.
But in that case, $z_{Q+1}=z_i$ for exactly one value of $i\in \{1,\ldots,Q\}$, so that the first factor becomes
$$\alpha+(Q-1)Y-(Q-1)\beta-\alpha.$$ 
Thus in addition setting $Y=\beta$ ensures each summand is zero, so $V_Q(\alpha,\beta)=0$.\\

Finally, the coefficient of $Y$ in $V_Q(\alpha,Y)$ is the product of $Q-1$ and
\begin{multline*}
\EE\left( \prod_{1\le i< j\le Q}(\alpha-X_{ij}) \right )\\
=\sum_{z\in \{1,\dots,Q\}^{Q}} \left (\prod_{1\le k\le Q} \pi_{z_k} \right ) \prod_{1\le i< j\le Q}\EE(\alpha-X_{ij}~|~Z_i=z_i,Z_j=z_j)  .
\end{multline*}
But  $\prod_{1\le i< j\le Q}\EE(\alpha-X_{ij}~|~ Z_i=z_i,Z_j=z_j)$ vanishes for all $z$
except  possibly for  those in  which all $z_i$, $1\le i\le Q$, are distinct, in which case it takes the value $(\alpha-\beta)^{\binom Q2}$. 
So the coefficent becomes
$$(Q-1)(Q!) \left (\prod_{1\le k\le Q} \pi_k \right ) (\alpha-\beta)^{\binom{Q}{2} }.$$
This is zero if, and only if, $\alpha=\beta$.
\end{proof}

\ 

%%%%%%%

\begin{proof}[Proof of  Theorem~\ref{thm:affil_Q=2}]

Since $\alpha$ is a real root of the cubic polynomial $U_2(X)$, to show $\alpha$ is uniquely identifiable it is enough to show that $\frac d{dX}U_2(X) \ge 0$. But
$$\frac d{dX}U_2(X)=3X^2-6m_1X+3m_2=3\left ( (X^2-m_1)^2+(m_2-m_1^2)\right ).$$
But $m_2-m_1^2\ge0$ because, using the Cauchy-Schwarz inequality,  
\begin{multline*}
m_2=\mathbb{E} (X_{ij}X_{ik}) =\mathbb{E}[\mathbb{E}(X_{ij}|Z_i)
\mathbb{E}(X_{ik}|Z_i) ] \\
= \mathbb{E}[\mathbb{E}(X_{ij}|Z_i)^2] \geq 
[\mathbb{E}(\mathbb{E}(X_{ij}|Z_i))]^2 =m_1^2.
\end{multline*} 

With $\alpha$ identified, since $\alpha\ne\beta$, we may uniquely recover $\beta$ as the root of the linear polynomial $V_2(\alpha,Y)$ with nonzero leading coefficient.

\end{proof}

\begin{proof}[Proof of Theorem~\ref{thm:poly_ST}]
Using equation \eqref{eq:m_1} to eliminate $\alpha$ from equations \eqref{eq:m_3} and \eqref{eq:m_2} respectively, gives two equations
\begin{equation*}
\begin{array}{ccl}
R(\beta) &=& a\beta^3+b\beta^2+c\beta +d=0,\\
S(\beta) &=& A\beta^2+B\beta+C =0,
  \end{array}
 \end{equation*} 
where 
\begin{equation*}
\left\{
  \begin{array}{ccl}
  a&=& -2s_2^3  + 3s_2s_3 -s_3, \\
b &=& {3m_1(s_2^3-2s_2s_3+s_3)},\\
c &=& 3m_1^2s_3(s_2-1),\\
d &=& m_1^3s_3-m_3s_2^3,\\   
  \end{array}
\right. 
\text{and}
\left\{
  \begin{array}{ccl}
A &=& s_3-s_2^2,\\
B &=& -2m_1(s_3-s_2^2),\\
C &=& m_1^2s_3-m_2s_2^2. 
  \end{array}
\right.   
\end{equation*}

To understand the degrees of these polynomials we need the following.

\begin{lem} \label{lem:aA}
Suppose $\boldsymbol \pi\in [0,1]^Q$ with $\sum_{q=1}^Q \pi_q=1$.   
\begin{itemize} 
\item[i)]  If $\pi_q>0$ for at least two values of $q$, then $a\ne 0$. 
\item[ii)] $A=0$ if, and only if, $\boldsymbol \pi$ is uniform on its support.
  \end{itemize}

\end{lem}

\begin{proof}
  To establish claim $i)$,  first observe that $0<s_2<  1$. Moreover, since
$s_3^2\le s_2s_4$ by the Cauchy-Schwarz inequality, and $s_4<s_2^2$ by comparing terms (since at least two $\pi_q>0$), we have $s_3<s_2^{3/2}$.
If $-2 s_2^3+3 s_2 s_3-s_3=0$, then $$s_2^{3/2}>s_3=\frac {2s_2^3}{3s_2-1},$$ where the denominator must be positive. Thus
$$1>\frac {2s_2^{3/2}}{3s_2-1},$$ so
$$0> {2s_2^{3/2}}- 3s_2+1.$$
However, the function $x\mapsto 2x^{3/2}- 3x+1$ is positive on $(0,1)$, so this is a contradiction.

Turning to claim  $ii)$,  we  have  $A=s_3-s_2^2$ and  by  the
Cauchy-Schwarz inequality,  $s_2^2=(\sum_q \pi_q^{3/2}\pi_q^{1/2})^2\le s_3,$
with      equality      if,      and      only     if, 
$(\pi_1^{3/2},\ldots,\pi_Q^{3/2})=\lambda (\pi_1^{1/2},\ldots, \pi_Q^{1/2})$ for some value $\lambda\in
\mathbb R$. This can only occur if on its support $\boldsymbol \pi$ is uniform.
\end{proof}
 
Returning to the proof of Theorem~\ref{thm:poly_ST}, if 
$\boldsymbol \pi$ is not uniform, we thus have $A\ne 0$ and dividing the polynomial $R(\beta)$ by $S(\beta)$ produces a linear remainder $T(\beta)$, which is calculated to be
\begin{multline*}
T(\beta)=\frac {s_2^2}{s_2^2-s_3}\left [ (m_2-m_1^2)(s_3-3s_3s_2+2s_2^3)\beta\right .\\
 \left .+   (s_3-s_2s_3)m_1^3 +(s_2^3-s_3)m_2m_1+(s_3s_2-s_2^3)m_3 \right ].
 \end{multline*}
Since any common zero of $R(\beta)$ and $S(\beta)$ must also be a zero of  $T(\beta)$, we can recover the parameters $\beta$ and $\alpha$ via the rational formulas
\begin{align}
\beta &=\frac{(s_3-s_2s_3)m_1^3 +(s_2^3-s_3)m_2m_1+(s_3s_2-s_2^3)m_3 }{(m_1^2-m_2)(2s_2^3-3s_3s_2+s_3)},\label{eq:beta}\\
\alpha &=\frac {m_1+(s_2-1)\beta}{s_2}.\label{eq:alpha}\end{align}
Note that a calculation shows 
\begin{equation}m_1^2-m_2=(\alpha-\beta)^2(s_2^2-s_3),\label{eq:m1m2}\end{equation}
which, since $A\ne 0$, is only zero in the trivial case of $\alpha=\beta$. Otherwise, since  $2s_2^3-3s_3s_2+s_3=-a\ne 0$ by part $i)$ of Lemma \ref{lem:aA}, the formulas \eqref{eq:beta} and \eqref{eq:alpha} are valid.

Equation \eqref{eq:m1m2}, together with part $ii)$ of Lemma \ref{lem:aA} further shows that if $m_2\ne m_1^2$, then $\boldsymbol \pi$ is not uniform.

\smallskip

If $m_2= m_1^2$, then $\boldsymbol \pi$ is uniform, and $S(\beta)$ is identically zero. However, in this case the coefficients of $$\tilde R(\beta)=\frac {Q^3}{1-Q}R(\beta)=\beta^3+\tilde b\beta^2+\tilde c \beta+\tilde d$$ simplify to
$$ \tilde b=-3m_1, \ \ \ \ \ \  \tilde c=3m_1^2,$$ 
$$\tilde d=\frac{Qm_1^3-m_3}{1-Q}=-m_1^3+       \frac {m_1^3-m_3}{1-Q}.$$
Thus
$$ \tilde R(\beta)=(\beta-m_1)^3+ \frac {m_1^3-m_3}{1-Q},$$
which has a unique real root
$$\beta=m_1+\left ( \frac {m_1^3-m_3}{Q-1}  \right )^{1/3}.$$
The parameter $\alpha$ can then be found by formula \eqref{eq:alpha}.

\end{proof}

%%%%%%%

\begin{proof}[Proof of Proposition~\ref{prop:bound_on_m}]

First, note that the distribution of $K_n$ may be parameterized using the 
elementary symmetric polynomials $\sigma_i$ evaluated at the $\{\pi_q\}_{1\le q \le Q}$, 
instead of the values $\{\pi_q\}_{1\le q \le Q}$. Indeed, the affiliation model 
distribution only involves the  $\pi_q$s through the symmetric expressions 
$$\sum_{q_1,\ldots,q_s,\atop q_i\ne q_j}\pi_{q_1}^{i_1}\ldots\pi_{q_s}^{i_s},$$ 
with $s\le Q$ and $\sum_{k\le s} i_k=n$, and these sums may be expressed 
as polynomials in the $\{\sigma_i(\pi_1,\dots,\pi_Q)\}_{1\le i\le n}$.  Thus for 
identifiability of the $\{\pi_q\}$ from the distribution of $K_n$, it is necessary 
that the $\{\pi_q\}$ be identifiable from the $\{\sigma_i(\pi_1,\dots,\pi_Q)\}_{1\le i\le n}$. 
Note also that $\sigma_1(\pi_1,\dots,\pi_Q)=\sum_{q=1}^Q \pi_i=1$ carries no information on the $\pi_q$s that is not already known.

Now if $n<Q$, identifying  $Q-1$ independent choices of the $\pi_q$ from the values of $n-1$ continuous functions of those $\pi_q$
is impossible.
\end{proof}

\begin{lem}\label{lem:4m} For the random graph affiliation model on $Q$ nodes, with binary edge state variables, uniform group priors, and connectivities $\alpha\ne \beta$, the moment inequality
$m_{41}>m_1^4$ holds.
\end{lem}
\begin{proof} Note
\begin{multline*}
m_{41}= \mathbb{E} [ \mathbb{E}(X_{12}X_{23} | Z_{1},Z_{3})
\mathbb{E}(X_{34}X_{41} | Z_{1},Z_{3}) ]
= \mathbb{E} [ \mathbb{E}(X_{12}X_{23} | Z_{1},Z_{3})^2]\\
\geq (\mathbb{E}[ \mathbb{E}(X_{12}X_{23} | Z_{1},Z_{3})])^2
= m_2^2. 
\end{multline*}
However,
equality occurs above only if $\mathbb{E}(X_{12}X_{23} | Z_{1},Z_{3})$ is constant. But
$$\mathbb{E}(X_{12}X_{23} | Z_{1}=i=Z_{3})=\frac{1}{Q}{\alpha^2}+\frac{Q-1}{Q}\beta^2,$$
$$\mathbb{E}(X_{12}X_{23} | Z_{1}=i\ne j=Z_{3})=\frac{2}{Q}{\alpha\beta}+\frac{Q-2}{Q}\beta^2,$$
so the difference of these expectations is $(\alpha-\beta)^2/Q\ne 0$.
Thus $m_{41}> m_2^2$.

A similar  argument that $m_{2} \ge  m_1^2$ was given in  the proof of
Theorem \ref{thm:affil_Q=2}, so the claim is established.
\end{proof}

%%%%%%%
\subsection{Proofs for the continuous parametric model}

\begin{proof}[Proof of Theorem~\ref{thm:parametric}]
With $\bar p_{q\ell}=1-p_{q\ell}$, the distribution of $(X_{ij},X_{ik},X_{jk})$ is given by
the mixture
\begin{multline}
  \sum_{1\le         q,\ell,m\le         Q}         \pi_q\pi_\ell\pi_m
  [\bar p_{q\ell}\delta_0(X_{ij})+p_{q\ell} F(X_{ij},\theta_{q\ell})] \times
[\bar p_{qm}\delta_0(X_{ik})+p_{qm} F(X_{ik},\theta_{qm})]\\
\times[\bar p_{\ell m}\delta_0(X_{jk})+p_{\ell m} F(X_{jk},\theta_{\ell m})].\label{eq:3dist}
\end{multline}

Since the distributions
$F(\cdot,\theta)$ have no point masses at $0$ by Assumption~\ref{hyp:ident}, the family $\mathcal M\cup \{\delta_0\}$
has identifiable parameters for finite mixtures, so Theorem 1 of \citet{Teicher67} applies to it.
Thus multiplying out the terms of the mixture in \eqref{eq:3dist} to view it as a mixture of products from $\mathcal M\cup \{\delta_0\}$, and
noting that by  Assumption \ref{hyp:theta_neq} certain of the components arise from unique choices of $q,\ell,m$
we can identify the terms of the form 
\begin{equation*}\label{eq:mixture_vectors}
     \pi_q\pi_\ell\pi_m   p_{q\ell}p_{qm}
   p_{\ell m} F(X_{ij},\theta_{q\ell}) F(X_{ik},\theta_{qm}) F(X_{jk},\theta_{\ell m}),
\end{equation*}
and the vectors in 
$$
\mathcal{C}=\{(\pi_q\pi_\ell\pi_mp_{q\ell}p_{qm}
   p_{\ell m} ; \theta_{q\ell},\theta_{qm},\theta_{\ell m}) ~|~ 1\le q,\ell,m\le Q\},
$$ 
but only as an unordered set. 
But by  Assumption~\ref{hyp:theta_neq}, there are only $Q$ vectors in
this       set       for       which      the       last       entries
$(\theta_{q\ell},\theta_{qm},\theta_{\ell m})$  are all equal.  Indeed, these
entries  are of  the form  $(\theta_{qq},\theta_{qq},\theta_{qq})$ for
some $1\le q\le Q$,  since the case where  these entries would  be of
the  form  $(\theta_{q\ell},\theta_{q\ell},\theta_{q\ell})$  for  some
$q\neq \ell$ is not possible.   Thus the $\theta_{qq}$ for $1\le q\le Q$ may be identified
as   well   as  the   corresponding   weights  $(\pi_qp_{qq})^3$,   or
equivalently the values $\pi_qp_{qq}$. 

Now, among the vectors in $\mathcal{C}$, exactly $3Q(Q-1)$ of them have two
of the last three entries equal. These entries are, up to order, of the
form   $(\theta_{qq},\theta_{q\ell},\theta_{q\ell})$,   for any  $q\neq  \ell$. 
Thus we obtain the set $\{(\pi_q^2\pi_\ell p_{q\ell}^2p_{qq}; \theta_{qq},\theta_{q\ell},\theta_{q\ell})\}_{1\le
  q< \ell \le Q}$, without regard to order.
Since we already identified the pairs $(\pi_qp_{qq}, \theta_{qq})$, we
may  take   the    ratio   between   the   weights   $\pi_q^2\pi_\ell
p_{q\ell}^2p_{qq}$   and   $\pi_qp_{qq}$   to   recover   the   values
$\pi_q\pi_\ell   p_{q\ell}^2$.    Thus   we   identify    the   set
$\{ (\pi_q\pi_\ell p_{q\ell}^2;\theta_{qq},\theta_{q\ell},\theta_{q\ell})\}_ {1\le
  q< \ell \le Q}$. 

Among these vectors,  we can match the ones whose  two last entries are 
equal, namely those of the form
$(\pi_q\pi_\ell p_{q\ell}^2; \theta_{qq}, \theta_{q\ell},\theta_{q\ell})$ with
$(\pi_q\pi_\ell p_{q\ell}^2; \theta_{\ell\ell},\theta_{q\ell},\theta_{q\ell})$.  This enables us
to  recover  the  values  $\theta_{q\ell}$, for $1\le
  q,\ell\le Q$.

By marginalizing the distribution of $(X_{ij},X_{ik},X_{jk})$, we also have the distribution of a single edge variable $X_{ij}$,
\begin{equation}
  \sum_{1\le         q,\ell\le         Q}         \pi_q\pi_\ell 
   [\bar p_{q\ell}\delta_0(X_{ij})+p_{q\ell} F(X_{ij},\theta_{q\ell})] 
  .\label{eq:1dist}
\end{equation}
and thus by our hypotheses can also identify $\{(\pi_q\pi_\ell p_{q\ell}, \theta_{q\ell})\}_{1\le q\le \ell\le Q}$, without order. But as the $\theta_{q\ell}$ have already been identified, we may use this to match $\pi_q\pi_\ell p_{q\ell}$ with $\pi_q\pi_\ell p_{q\ell}^2$ and thus recover $p_{q\ell}$ from the ratio. From $\pi_qp_{qq}$ and $p_{qq}$ we can then recover $\pi_q$.

%,as well as $\pi_q\pi_\ell$.

%Returning to the set $\mathcal C$, from the vectors whose last three entries are all different, we may recover the weights
%$\pi_q\pi_\ell   \pi_m   p_{q\ell}   p_{qm}  p_{\ell   m}$ for $1\le q<\ell<m\le Q$. But since $p_{q\ell} ,  p_{qm},  p_{\ell   m}$ and 
%$\pi_q\pi_\ell, \pi_q\pi_m, \pi_\ell\pi_m$ are already known, this gives $\pi_i$ for $1\le i\le Q$.

Thus, all parameters of the model are identified, up to permutation on the group labels.
\end{proof}

\begin{proof}[Proof of Theorem~\ref{thm:parametric_affiliation}]
From   the  distribution  of  $K_3$,  we  can
distinguish          $(\alpha,         \theta_\text{in})$         from
$(\beta,\theta_{\text{out}})$ as follows: The distribution of $K_3$ is
the mixture of either $4$ (when $Q=2$) or $5$ (when $Q\ge 3$) different 3-dimensional components. Since the distributions
$F(\cdot,\theta)$ do not have point masses at $0$ by Assumption~\ref{hyp:ident}, we can identify
from this  mixture that part with no such Dirac masses  in it, which is the mixture
\begin{multline*}
\alpha^3   \Big(\sum_{q=1}^Q\pi_q^3\Big)F(\cdot,\theta_{\text {in}  })
\otimes  F(\cdot,\theta_{\text {in} })  \otimes  F(\cdot,\theta_{\text
 {in}}) \\
 +\alpha\beta^2    \Big(\sum_{1\le   q\neq    \ell\le   Q}
\pi_q^2\pi_\ell     \Big)F(\cdot,\theta_{\text {in}    })    \otimes
F(\cdot,\theta_{\text {out} }) \otimes F(\cdot,\theta_{\text {out} }) \\
+\alpha\beta^2    \Big(\sum_{1\le   q\neq    \ell\le   Q}
\pi_q^2\pi_\ell     \Big)F(\cdot,\theta_{\text  {out}    })    \otimes
F(\cdot,\theta_{\text {in} }) \otimes F(\cdot,\theta_{\text {out} }) \\
+\alpha\beta^2    \Big(\sum_{1\le   q\neq    \ell\le   Q}
\pi_q^2\pi_\ell     \Big)F(\cdot,\theta_{\text {out}    })    \otimes
F(\cdot,\theta_{\text {out} }) \otimes F(\cdot,\theta_{\text{in} }) \\
+\beta^3 \Big(\sum_{ q,\ell,m \text{ distinct}}
\pi_q\pi_\ell\pi_m     \Big)F(\cdot,\theta_{\text{out}    })    \otimes
F(\cdot,\theta_{\text{out} }) \otimes F(\cdot,\theta_{\text{out} }),
\end{multline*}
where the last term appears only when $Q\ge 3$.

By  Theorem 1 of \citet{Teicher67}  and Assumption~\ref{hyp:ident},  this 3-dimensional mixture  has identifiable
parameters, up to label swapping issues.  At most two terms in this
mixture  have the  same  measure  $F$ in  each  coordinate. The  three
remaining  terms  have  two   coordinates  which  are  equal, involving $\theta_{\text{out}}$,  and  one
different, involving $\theta_{\text{in}}$. Thus we can distinguish between $\theta_{\text{in}}$
and $\theta_{\text{out}}$.

We may also determine $\alpha^3(\sum_q \pi_q^3)$ as the weight of $F(\cdot,\theta_{\text{in}})\otimes
F(\cdot,\theta_{\text{in}})\otimes F(\cdot,\theta_{\text{in}})$.
Similarly from the $\delta_0 \otimes F(\cdot,\theta_{\text{in}}) \otimes F(\cdot,\theta_{\text{in}})$ term in the full mixture, we may
recover the weight $(1-\alpha)\alpha^2 (\sum_q \pi_q^3)$.
Summing these two weights yields $\alpha^2 (\sum_q \pi_q^3)$,
and then dividing the first by this, we recover $\alpha$.

The parameter $\beta$ is similarly recovered from the weights of $F(\cdot,\theta_{\text{out}})\otimes
F(\cdot,\theta_{\text{out}})\otimes F(\cdot,\theta_{\text{in}})$ and $\delta_0\otimes
F(\cdot,\theta_{\text{out}})\otimes F(\cdot,\theta_{\text{in}})$.

Next we consider the distribution of $K_n$ for various $n$.  This is a mixture of
many different $\binom  {n}{2}$-dimensional components. As above, we can identify up to label swapping the components with no $\delta_0$ factors
in   this   mixture.  But   as   we   already   know  the   value   of
$\theta_{\text{in}}$, we  can identify the  term $\otimes_{1\le i<j\le
  n}  F(X_{ij},\theta_{\text{in}})$  in  this  mixture, and  thus  its
corresponding prior $\alpha^n  \sum_{q} \pi_{q}^n$. Since $\alpha$ has
been previously identified, this uniquely determines $\sum_{q}
\pi_{q}^n$. Note that using the distribution of
$K_Q$, we can obtain the distribution  of each $K_n$ with $n\le Q$ and
thus the values $\{\sum_{q} \pi_{q}^n\}_{n\le Q}$. 

By the Newton identities, these values determine the values of
elementary symmetric polynomials $\{\sigma_{n}( \pi_1,\dots,\pi_Q)\}_{n\le Q}$. These, in turn, are (up to sign) the coefficients of the monic polynomial whose roots (with multiplicities) are precisely
$\{\pi_q\}_{ 1\le q\le Q}$.  Thus the node priors are determined, up to order.

\end{proof}

%%%%%%%

\subsection{Proof of Theorem~\ref{thm:discrete}}

The proof follows the strategy described in Section~\ref{sec:scheme}.
We thus proceed with a base case, an extension step, and a conclusion. 

\paragraph*{Base case}
We consider a subset $\mathcal{E}$ of the set of all edges
over  $m$ vertices, with $m$ and $\mathcal{E}$ to be chosen later. Let $A$ be
the $Q^m\times \kappa^{|\mathcal{E}|}$ matrix containing the
  probabilities of the clumped random variable $Y=(X_e)_{ e \in
  \mathcal{E}}$ with state space $\{1,\ldots,\kappa\} ^{\vert \mathcal{E} \vert}$, conditional on the hidden states of the $m$ vertices.

Let $\mathcal{I}\in\{1,\ldots,Q\}^m$ be a vector specifying particular states of all the node variables. For each
edge $e \in \mathcal{E}$, the endpoints are in some set of
hidden states $\{q,l\}$, which we denote by $\mathcal{I}(e)$.
The $(\mathcal{I},(x_e)_{e\in \mathcal{E}})$-entry of the matrix $A$
is then given by  
\begin{equation*}
  \prod_{e\in \mathcal{E}}\prod_{k=1}^{\kappa} (p_{\mathcal{I}(e)}(k))^{1_{x_e=k}},
\end{equation*}
where $1_A$ is the indicator function for a set $A$.

For each edge $e$ in the graph, we introduce $\kappa$
indeterminates, $t_{e,1},\ldots,t_{e,\kappa}$. We create a
$\kappa^{|\mathcal{E}|}$-element column vector $\boldsymbol{t}$
indexed by the states of the clumped variable $Y$, whose
$(x_e)_{e\in \mathcal{E}}$-th entry is  given by 
\begin{equation*}
  \prod_{e\in \mathcal{E}} \prod_{k=1}^{\kappa} t_{e,k}^{1_{x_e=k}} . 
\end{equation*}
Then the $\mathcal{I}$th entry of the $Q^m$-entry vector
$A\boldsymbol{t}$ is the polynomial function 
\begin{equation*}
  f_{\mathcal{I}} = \sum_{(x_e)_{e\in \mathcal{E}}} \prod_{e \in
    \mathcal{E}} \prod_{k=1}^{\kappa}  \{p_{\mathcal{I}(e)}(k) t_{e,k}\}^{1_{x_e=k}}
= \prod_{e\in \mathcal{E}}
  \Big(p_{\mathcal{I}(e)}(1)t_{e,1} + \cdots +
  p_{\mathcal{I}(e)}(\kappa)t_{e,\kappa}\Big).
\end{equation*}

Independence of the rows of $A$ is equivalent to the independence of
the polynomials $\{f_{\mathcal{I}}\}_{\mathcal{I} \in \{1,\ldots, Q\}^m}$.
Thus, suppose that we have
\begin{equation}
  \sum_{\mathcal{I}} a_{\mathcal{I}}f_{\mathcal{I}} \equiv 0, \label{eq:depend}
\end{equation}
and let us show then that every $a_{\mathcal I}$ must be $0$.

For a specific $e\in \mathcal{E}$, and any choice $\{q,l\}$ with $1\le
q\le l\le Q$, one can choose a point 
$\mathbf{t}_{e,\{q,l\}}=(t_{e,1},\ldots,t_{e,\kappa}) \in \mathbb{R}^\kappa$
in the zero set of all the polynomial functions $f_{\mathcal{I}}$ in \eqref{eq:depend},
except those with $\mathcal{I}(e)=\{q,l\}$. 
To see this, let $M$ be the $\binom {Q+1}{2}\times \kappa$ matrix whose $\{q,l\}$th 
row is given by the vector $\mathbf{p}_{ql}=(p_{ql}(1),\ldots,p_{ql}(\kappa))$. 
$M$ has full row rank since its rows are independent by assumption. Thus 
there is a solution $\mathbf{t}_{e,\{q,l\}}$ to
\begin{equation*}
  M\mathbf{t}_{e,\{q,l\}} = \mathbf e_{\{q,l\}},
\end{equation*}
where $\mathbf e_{\{q,l\}}$ is the
vector of size $\binom {Q+1}{2}$ with zero entries, except the
$\{q,l\}$th which is equal to 1.  The independence assumption  also implies $\kappa \ge
  \binom{Q+1}{2}$.

  Note  that in this construction we have only specified group assignments to two nodes up to node permutation. Thus if the $\{q,l\}$ row of $M$ is related to an edge $e=(i,j)$ because $\mathcal{I}(e)=\{q,l\}$, we may have that either $i$ is in state $q$ and $j$ is in state $l$, or  $i$ is in state $l$ and $j$ is in state $q$.

By evaluating  the $f_{\mathcal I}$ at $\mathbf{t}_{e,\{q,l\}}$ for many edges $e$ and choices of node
states $\{q,l\}$, we can annihilate all the polynomials
$f_{\mathcal{I}}$ except those satisfying specific constraints on the node
states. More precisely, we can make vanish all the $f_{\mathcal{I}}$
except those for which $\mathcal{I}$ satisfies the condition that for some subset of edges  $\mathcal{E}'\subseteq\mathcal{E}$ and some sequence of unordered node assignments $(\{q_e,l_e\})_{e \in \mathcal{E}'}$ we have
\begin{equation}\label{eq:sets}
\mathcal{I} \in  \bigcap_{e \in \mathcal{E}'}
 \mathcal{S}(e;\{q_e,l_e\}) ,
\end{equation}
 where $\mathcal{S}\left(e;\{q_e,l_e\}\right)= 
\left\{ \mathcal{I}\in \{1,\ldots,Q\}^m \mid \mathcal{I}(e) =\{q_e,l_e\} \right\} $.

To conclude that each
$a_{\mathcal{I}}=0$ in equation \eqref{eq:depend}, it is enough to construct for 
every $\mathcal{I}\in \{1,\ldots,Q\}^m$ a set as in
\eqref{eq:sets} containing only  $\mathcal{I}$. 

In fact, this can be achieved with only $m=3$
vertices and the full set of edges
$\mathcal{E}=\{(1,2),(1,3),(2,3)\}$. Indeed, up to permutation of the nodes and of the
labels of the groups, $\mathcal{I}$ can take only three different values, namely $(1,1,1), (1,1,2)$ and $(1,2,3)$. Using a node assignment on the edges in $\mathcal{E}'=\{(1,2),(2,3)\}$, we get  
\begin{eqnarray*}
  \{(1,1,1)\} &=& \mathcal{S}\left((1,2);\{1,1\}) \cap  \mathcal{S}((2,3);\{1,1\}\right) \\
  \{(1,1,2)\} &=& \mathcal{S}\left((1,2);\{1,1\}) \cap  \mathcal{S}((2,3);\{1,2\}\right) \\
  \{(1,2,3)\} &=& \mathcal{S}\left((1,2);\{1,2\}) \cap  \mathcal{S}((2,3);\{2,3\}\right) .
\end{eqnarray*}

Thus, we proved the following lemma.
\begin{lem}\label{lem:base}
With  $\mathcal{E}$ the complete set of edges over
$m=3$ vertices,  the $Q^3\times \kappa^3$ matrix $A$ containing
the probabilities of the clumped variable $Y=(X_e)_{e \in
  \mathcal{E}}$, conditional on the hidden states $Z=(Z_1,Z_2,Z_3)
\in \{1,\ldots,Q\}^3$ has full row rank $Q^3$, provided the $\kappa$-entry vectors  $\{\mathbf{p}_{ql}\}_{1\le q\le l\le Q}$ are linearly independent. 
\end{lem}

\paragraph*{Conclusion of the proof}
The Lemma provides the base case, with the extension step of Section \ref{sec:scheme} then applying. Thus with $n=m^2=9$ nodes, Kruskal's Theorem may be
applied to identify, up to simultaneous row permutation, $\mathbf v$, $M_1$, $M_2$, and $M_3$ as defined in that section.

The rest of the proof follows the same lines as the conclusion in the proof of Theorem~\ref{thm:binary_Qge3_nonaff}, replacing the numbers $p_{ql}$ by the vectors $\mathbf{p}_{ql}$ and noting that these vectors are assumed to be linearly independent.

%%%%%%%%%%
\subsection{Proof of Theorem~\ref{thm:nonparametric}}

For convenience, we present the argument  assuming the state space of the $\mu_{ql}$ is a subset of $\mathbb R$. The more general situation of a multidimensional state space can be handled similarly, along the lines of the proof of Theorem 9 of  \cite{ECJ}.

Let $M_{ql}$ denote the c.d.f.~of $\mu_{ql}=(1-p_{ql})\delta_0+p_{ql}F_{ql}$. Since the measures $\{\mu_{ql} ~|~1\le q\le l\le Q\}$ are assumed to be linearly independent, so are
the functions $\{ M_{ql} ~|~ 1\le q\le l\le Q\}$.
Applying Lemma 17 of \citet{ECJ} to this set of functions, there exists some $\kappa\in \NN$ and cutpoints $u_1<u_2<\cdots<u_{\kappa-1}$ such that the vectors
$$\{ ( M_{ql}(u_1),M_{ql}(u_2),\dots,M_{ql}(u_{\kappa-1}),1)~|~1\le q\le l\le Q\}$$
are independent.  Note $\kappa\ge \binom {Q+1}2$. Also
by adding additional cutpoints if necessary, and thereby increasing $\kappa$, we may assume that among the $u_i$ are any specific real numbers we like.

The independence of the above vectors is equivalent to the independence of the vectors $\{\bar M_{ql} ~|~1\le q\le l\le Q\}$, where
$$\bar M_{ql}= \left( M_{ql}(u_1),M_{ql}(u_2)-M_{ql}(u_1),\dots,M_{ql}(u_{\kappa-1})-M_{ql}(u_{\kappa-2}),1-M_{ql}(u_{\kappa-1})\right).$$
Note that the $k$th entry of $\bar M_{ql}$ is simply the probability that a variable with distribution $\mu_{ql}$ takes values in the intervals $I_k=(u_{k-1},u_k]$ (with the convention that $u_{0}=-\infty, u_\kappa=\infty$).
To formalize this, let 
$$Y_{ij} =\sum_{k=1}^\kappa k 1_{I_k}(X_{ij})$$
be the random variable with state space $\{1, 2, \ldots, \kappa\}$
indicating the interval in which the value of $X_{ij}$ lies. Thus,
conditional on $Z_i=q,Z_j=l$, the random variables $X_{ij}$ and $Y_{ij}$
have respective c.d.f.s $M_{ql}$ and $\bar M_{ql}$.

Now from the distribution of the continuous random graph mixture model on $K_9$, with edge variables $(X_{ij})_{1\le i<j\le 9}$, by binning the values of the 36 edge variables into sets
of the form $\prod_{1\le i<j\le 9} I_{k_{ij}}$ with $1\le k_{ij}\le \kappa$, we obtain the distribution for the discrete edge variables $(Y_{ij})_{1\le i<j\le 9}$ of a random graph mixture model with the same group priors on the nodes, and with mixture components built from the distributions $\bar M_{ql}$ associated to $\mu_{ql}$. By Theorem \ref{thm:discrete}, the parameters of the discrete model are identifiable, up to label swapping. Imposing an arbitrary labeling, we have identified the node group priors $\pi_q$, $1\le q\le Q$, and for each pair of groups $q\le l$ the vector $\bar M_{ql}$. By summing entries of $\bar M_{ql}$, we obtain values of $M_{ql}(u_k)$ for $k=1,2,\dots, \kappa-1$. Since we may additionally determine
$M_{ql}(t)$ for any real number $t$ by including it as a cutpoint,  $M_{ql}$, and hence $\mu_{ql}$, is uniquely determined.

%%%%%%%%%%%%%%%%%%%%%%%%%%
\section{Acknowledgements}

The authors thank the Statistical and Applied Mathematical Sciences Institute for their support during residencies in which some of this work was undertaken. ESA and JAR also thank the Laboratoire Statistique et G\'enome for its hospitality.
JAR additionally thanks Universit\'e d'\'Evry Val d'Essonne for a Visiting Professorship during which this work was completed. 
ESA and JAR received support from the National Science Foundation, grant DMS 0714830, while CM has been supported  by the French Agence Nationale de la Recherche under grant NeMo ANR-08-BLAN-0304-01.

\bibliographystyle{elsarticle-harv}
\bibliography{mixnet_ident}

\begin{thebibliography}{36}
\expandafter\ifx\csname natexlab\endcsname\relax\def\natexlab#1{#1}\fi
\expandafter\ifx\csname url\endcsname\relax
  \def\url#1{\texttt{#1}}\fi
\expandafter\ifx\csname urlprefix\endcsname\relax\def\urlprefix{URL }\fi

\bibitem[{Airoldi et~al.(2008)Airoldi, Blei, Fienberg, and Xing}]{Airoldi}
Airoldi, E., Blei, D., Fienberg, S., Xing, E., 2008. Mixed-membership
  stochastic blockmodels. Journal of Machine Learning Research 9, 1981--2014.

\bibitem[{Allman et~al.(2009)Allman, Matias, and Rhodes}]{ECJ}
Allman, E., Matias, C., Rhodes, J., 2009. Identifiability of parameters in
  latent structure models with many observed variables. Ann. Statist. 37~(6A),
  3099--3132.

\bibitem[{Ambroise and Matias(2010)}]{Matias_Ambroise}
Ambroise, C., Matias, C., 2010. New consistent and asymptotically normal
  estimators for random graph mixture models. Tech. rep., arXiv:1003.5165.

\bibitem[{Barrat et~al.(2004)Barrat, Barth{\'e}lemy, Pastor-Satorras, and
  Vespignani}]{PNAS_Barrat}
Barrat, A., Barth{\'e}lemy, M., Pastor-Satorras, R., Vespignani, A., 2004. {The
  architecture of complex weighted networks}. PNAS 101~(11), 3747--3752.

\bibitem[{Berge(1976)}]{Berge}
Berge, C., 1976. {Graphs and hypergraphs. Translated by Edward Minieka. 2nd
  rev. ed.} {North-Holland Mathematical Library. Vol. 6. Amsterdam - Oxford:
  North- Holland Publishing Company; New York:American Elsevier Publishing}.

\bibitem[{Carreira-Perpi{\~n}{\'a}n and Renals(2000)}]{Carreira}
Carreira-Perpi{\~n}{\'a}n, M., Renals, S., 2000. Practical identifiability of
  finite mixtures of multivariate {B}ernoulli distributions. Neural Comp.
  12~(1), 141--152.

\bibitem[{Cox et~al.(1997)Cox, Little, and O'Shea}]{CLO}
Cox, D., Little, J., O'Shea, D., 1997. Ideals, varieties, and algorithms, 2nd
  Edition. Springer-Verlag, New York.

\bibitem[{Daudin et~al.(2008)Daudin, Picard, and Robin}]{Daudin}
Daudin, J.-J., Picard, F., Robin, S., 2008. A mixture model for random graphs.
  Statist. Comput. 18~(2), 173--183.

\bibitem[{Daudin et~al.(2010)Daudin, Pierre, and Vacher}]{Daudin10}
Daudin, J.-J., Pierre, L., Vacher, C., 2010. Model for heterogeneous random
  networks using continuous latent variables and an application to a
  tree-fungus network. Biometrics, to appear.

\bibitem[{Erd\H{o}s and Gallai(1961)}]{Erdos_Gallai}
Erd\H{o}s, P., Gallai, T., 1961. {Graphs with points of prescribed degree.
  (Graphen mit Punkten vorgeschriebenen Grades.)}. Mat. Lapok 11, 264--274.

\bibitem[{Erd{\H{o}}s and R{\'e}nyi(1959)}]{ER1}
Erd{\H{o}}s, P., R{\'e}nyi, A., 1959. On random graphs. {I}. Publ. Math.
  Debrecen 6, 290--297.

\bibitem[{Frank and Harary(1982)}]{FH82}
Frank, O., Harary, F., 1982. Cluster inference by using transitivity indices in
  empirical graphs. J. Amer. Statist. Assoc. 77~(380), 835--840.

\bibitem[{Gyllenberg et~al.(1994)Gyllenberg, Koski, Reilink, and
  Verlaan}]{Gyllenberg}
Gyllenberg, M., Koski, T., Reilink, E., Verlaan, M., 1994. Nonuniqueness in
  probabilistic numerical identification of bacteria. J. Appl. Probab. 31~(2),
  542--548.

\bibitem[{Handcock et~al.(2007)Handcock, Raftery, and Tantrum}]{Handcock}
Handcock, M., Raftery, A., Tantrum, J., 2007. Model-based clustering for social
  networks. J. Roy. Statist. Soc. Ser. A 170~(2), 301--354.

\bibitem[{Holland et~al.(1983)Holland, Laskey, and Leinhardt}]{Holland_etal_83}
Holland, P., Laskey, K., Leinhardt, S., 1983. Stochastic blockmodels: some
  first steps. Social networks 5, 109--137.

\bibitem[{Kruskal(1976)}]{Kruskal76}
Kruskal, J., 1976. More factors than subjects, tests and treatments: an
  indeterminacy theorem for canonical decomposition and individual differences
  scaling. Psychometrika 41~(3), 281--293.

\bibitem[{Kruskal(1977)}]{Kruskal77}
Kruskal, J., 1977. Three-way arrays: rank and uniqueness of trilinear
  decompositions, with application to arithmetic complexity and statistics.
  Linear Algebra and Appl. 18~(2), 95--138.

\bibitem[{Latouche et~al.(2009)Latouche, Birmel{\'e}, and
  Ambroise}]{Latouche09}
Latouche, P., Birmel{\'e}, E., Ambroise, C., 2009. Overlapping stochastic block
  models. Tech. rep., arXiv:0910.2098.

\bibitem[{Mariadassou and Robin(2010)}]{Mariadassou_Robin}
Mariadassou, M., Robin, S., 2010. Uncovering latent structure in valued graphs:
  a variational approach. Annals of Applied Statistics, to appear.

\bibitem[{McLachlan and Peel(2000)}]{MacLachlan_Peel}
McLachlan, G., Peel, D., 2000. Finite mixture models. Wiley Series in
  Probability and Statistics: Applied Probability and Statistics.
  Wiley-Interscience, New York.

\bibitem[{Newman(2003)}]{Newman}
Newman, M. E.~J., 2003. The structure and function of complex networks. SIAM
  Rev. 45~(2), 167--256 (electronic).

\bibitem[{Newman(2004)}]{Newman_weighted}
Newman, M. E.~J., 2004. Analysis of weighted networks. Phys. Rev. E 70, 056131.

\bibitem[{Newman and Leicht(2007)}]{Newman_Leicht}
Newman, M. E.~J., Leicht, E.~A., 2007. {Mixture models and exploratory analysis
  in networks}. PNAS 104~(23), 9564--9569.

\bibitem[{Nowicki and Snijders(2001)}]{NS01}
Nowicki, K., Snijders, T., 2001. Estimation and prediction for stochastic
  blockstructures. J. Amer. Statist. Assoc. 96~(455), 1077--1087.

\bibitem[{Petrie(1969)}]{Petrie69}
Petrie, T., 1969. Probabilistic functions of finite state {M}arkov chains. Ann.
  Math. Statist 40, 97--115.

\bibitem[{Picard et~al.(2009)Picard, Miele, Daudin, Cottret, and
  Robin}]{Picard_BMC}
Picard, F., Miele, V., Daudin, J.-J., Cottret, L., Robin, S., 2009. Deciphering
  the connectivity structure of biological networks using {M}ix{N}et. BMC
  Bioinformatics 10, 1--11.

\bibitem[{Rhodes(2010)}]{John_on_Kruskal}
Rhodes, J., 2010. A concise proof of {K}ruskal's theorem on tensor
  decomposition. Linear Algebra and its Applications 432~(7), 1818--1824.

\bibitem[{Snijders and Nowicki(1997)}]{SN97}
Snijders, T., Nowicki, K., 1997. Estimation and prediction for stochastic
  blockmodels for graphs with latent block structure. J. Classification 14~(1),
  75--100.

\bibitem[{Tallberg(2005)}]{Tallberg05}
Tallberg, C., 2005. A {B}ayesian approach to modeling stochastic
  blockstructures with covariates. Journal of Mathematical Sociology 29~(1),
  1--23.

\bibitem[{Teicher(1961)}]{Teicher61}
Teicher, H., 1961. Identifiability of mixtures. Ann. Math. Statist. 32,
  244--248.

\bibitem[{Teicher(1963)}]{Teicher63}
Teicher, H., 1963. Identifiability of finite mixtures. Ann. Math. Statist. 34,
  1265--1269.

\bibitem[{Teicher(1967)}]{Teicher67}
Teicher, H., 1967. Identifiability of mixtures of product measures. Ann. Math.
  Statist 38, 1300--1302.

\bibitem[{Tomasi and Bro(2006)}]{Tomasi}
Tomasi, G., Bro, R., 2006. A comparison of algorithms for fitting the {PARAFAC}
  model. Comput. Statist. Data Anal. 50~(7), 1700--1734.

\bibitem[{White et~al.(1976)White, Boorman, and Breiger}]{White_76}
White, H., Boorman, S., Breiger, R., 1976. Social structure from multiple
  networks i: Blockmodels of roles and positions. American Journal of Sociology
  81, 730--779.

\bibitem[{Zanghi et~al.(2008)Zanghi, Ambroise, and Miele}]{Zanghi08}
Zanghi, H., Ambroise, C., Miele, V., 2008. Fast online graph clustering via
  {E}rd{\H{o}}s {R}{\'e}nyi mixture. Pattern Recognition 41~(12), 3592--3599.

\bibitem[{Zanghi et~al.(2010)Zanghi, Picard, Miele, and Ambroise}]{Zanghi09}
Zanghi, H., Picard, F., Miele, V., Ambroise, C., 2010. Strategies for online
  inference of network mixture. Annals of Applied Statistics, to appear.

\end{thebibliography}

\end{document}